\def\yes{\if00}
\def\no{\if01}
\def\iftwelvept{\yes}
\def\ifpsfont{\no}

\iftwelvept
\documentclass[leqno,12pt]{amsart}
\else
\documentclass[leqno]{amsart}
\fi
\usepackage{amssymb}
\usepackage{amscd}
\usepackage{latexsym}
\usepackage{verbatim}
\ifpsfont
\usepackage[T1]{fontenc}
\usepackage{times}
\else\fi

\iftwelvept
\else\fi

\theoremstyle{plain}
\newtheorem{Theorem}{Theorem}[section]
\newtheorem{Proposition}[Theorem]{Proposition}
\newtheorem{Lemma}[Theorem]{Lemma}
\newtheorem{Corollary}[Theorem]{Corollary}
\newtheorem{Claim}{Claim}[Theorem]

\theoremstyle{definition}
\newtheorem{Definition}[Theorem]{Definition}

\renewcommand{\theTheorem}{\arabic{section}.\arabic{Theorem}}
\renewcommand{\theClaim}{\arabic{section}.\arabic{Theorem}.\arabic{Claim}}
\renewcommand{\theequation}{\arabic{section}.\arabic{Theorem}.\arabic{Claim}}

\def\rom{\textup}
\newcommand{\ZZ}{{\mathbb{Z}}}
\newcommand{\QQ}{{\mathbb{Q}}}

\newcommand{\CC}{{\mathbb{C}}}

\newcommand{\PP}{{\mathbb{P}}}

\newcommand{\OO}{{\mathcal{O}}}
\newcommand{\XX}{{\mathcal{X}}}
\newcommand{\YY}{{\mathcal{Y}}}
\newcommand{\ZZZ}{{\mathcal{Z}}}
\newcommand{\LL}{{\mathcal{L}}}
\newcommand{\MM}{{\mathcal{M}}}
\newcommand{\codim}{\operatorname{codim}}

\newcommand{\Supp}{\operatorname{Supp}}

\newcommand{\Pic}{\operatorname{Pic}}

\newcommand{\zero}{\operatorname{div}}
\newcommand{\Proof}{{\sl Proof.}\quad}
\newcommand{\adeg}{\widehat{\operatorname{deg}}}
\newcommand{\trdeg}{\operatorname{tr.deg}}
\newcommand{\rank}{\operatorname{rk}}
\newcommand{\acherncl}{\widehat{{c}}}

\newcommand{\QED}{{\unskip\nobreak\hfil\penalty50\quad\null\nobreak\hfil
{$\Box$}\parfillskip0pt\finalhyphendemerits0\par\medskip}}
\newcommand{\rest}[2]{\left.{#1}\right\vert_{{#2}}}
\begin{document}

%%%%%%%%%%%
%% Title %%
%%%%%%%%%%%
\title[The canonical arithmetic height over a finitely generated field]%
{The canonical arithmetic height of subvarieties \\
of an abelian variety \\
over a finitely generated field}
\author{Atsushi Moriwaki}
\address{Department of Mathematics, Faculty of Science,
Kyoto University, Kyoto, 606-8502, Japan}
\email{moriwaki@kusm.kyoto-u.ac.jp}
\date{19/October/1999, 0:20PM (JP), (Version 1.0)}
\keywords{height function, abelian variety,
finitely generated field,
Arakelov Geometry, Bogomolov's conjecture}
\subjclass{Primary 11G35, 14G25, 14G40; Secondary 11G10, 14K15}
%\begin{abstract}
%\end{abstract}

%%%%%
% 1st Draft
%%%%%
%\vskip -1cm
%\hfill\fbox{{\large\bf 1st draft}} %
%\footnote{I welcome any comments and suggestions.}\par
%\vskip .5cm
%%%%%

\maketitle

%\tableofcontents

\section*{Introduction}
This paper is the sequel of \cite{MoArht}.
In \cite{ZhSmall}, S. Zhang defined the canonical
height of subvarieties of an abelian variety over a number field
in terms of adelic metrics.
In this paper, we generalize it to an abelian variety defined over
a finitely generated field over $\QQ$.
Our way is slightly different from his method.
Instead of using adelic metrics directly, we introduce an
adelic sequence and an adelic structure (cf. \S\S\ref{subsec:adelic:sequence}).

Let $K$ be a finitely generated field over $\QQ$
with $d = \trdeg_{\QQ}(K)$, 
and $\overline{B} = (B; \overline{H}_1, \ldots, \overline{H}_d)$
a polarization of $K$, i.e.,
$B$ is a projective arithmetic variety whose function field is $K$, and
$\overline{H}_1, \ldots, \overline{H}_d$ are 
nef $C^{\infty}$-hermitian line bundles on $B$.
Let $A$ be an abelian variety over $K$, and $L$ a symmetric ample line
bundle on $A$. Fix a projective arithmetic variety $\mathcal{A}$ over $B$ and
a nef $C^{\infty}$-hermitian $\QQ$-line bundle $\overline{\LL}$
on $\mathcal{A}$
such that $A$ is the generic fiber of $\mathcal{A} \to B$ and
$\LL$ is isomorphic to $L$ on $A$.
Then we can assign the naive height
$h^{\overline{B}}_{(\mathcal{A},\overline{\LL})}(X)$
to a subvariety $X$ of $A_{\overline{K}}$.
Indeed, if $X$ is defined over $K$,
$h^{\overline{B}}_{(\mathcal{A},\overline{\LL})}(X)$
is given by
\[
 \frac{\adeg\left( 
\acherncl_1(\rest{\overline{\LL}}{\XX})^{\cdot \dim X  +1} \cdot 
\acherncl_1(\pi_{\XX}^*(\overline{H}_1)) \cdots
\acherncl_1(\pi_{\XX}^*(\overline{H}_d))
\right)}{(\dim X + 1)\deg(\rest{L}{X}^{\dim X})},
\]
where $\XX$ is the Zariski closure of $X$ in $\mathcal{A}$ and
$\pi_{\XX} : \XX \to B$ is the canonical morphism.
The canonical height $\hat{h}_L^{\overline{B}}(X)$ of $X$ with respect to
$L$ and $\overline{B}$ is
characterized by the following properties:
\begin{enumerate}
\renewcommand{\labelenumi}{(\alph{enumi})}
\item
$\hat{h}_L^{\overline{B}}(X) \geq 0$
for all subvarieties $X$ of $A_{\overline{K}}$.

\item
There is a constant $C$ such that
\[
 \left\vert \hat{h}_L^{\overline{B}}(X) - 
        h^{\overline{B}}_{(\mathcal{A}, \overline{\LL})}(X) \right\vert \leq C
\]
for all subvarieties $X$ of $A_{\overline{K}}$.

\item
$\hat{h}_L^{\overline{B}}([N](X)) = N^2 \hat{h}_L^{\overline{B}}(X)$
for all subvarieties $X$ of $A_{\overline{K}}$ and
all non-zero integers $N$.
\end{enumerate}
The main result of this paper is the following theorem, which is
a generalization of \cite{ZhEqui}.

\renewcommand{\theTheorem}{\!\!}
\begin{Theorem}[cf. Theorem~\ref{thm:height:zero}]
If the polarization $\overline{B}$ is big \rom{(}i.e.,
$\overline{H}_1, \ldots, \overline{H}_d$ are nef and big\rom{)}, 
then, for a subvariety $X$ of $A_{\overline{K}}$,
the following are equivalent.
\begin{enumerate}
\renewcommand{\labelenumi}{(\arabic{enumi})}
\item
$X$ is a translation of an abelian subvariety by a torsion point.

\item
The set 
$\{ x \in X(\overline{K}) \mid \hat{h}_L^{\overline{B}}(x)
   \leq \epsilon \}$
is Zariski dense in $X$ for every $\epsilon > 0$.

\item
The canonical height of $X$ with respect to $L$ and $\overline{B}$
is zero, i.e.,
$\hat{h}^{\overline{B}}_{L}(X) = 0$.
\end{enumerate}
\end{Theorem}

Next let us consider a case where a curve and its Jacobian.
Let $X$ be a smooth projective curve of genus $g \geq 2$ over 
$K$, and $J$ the Jacobian of $X$.
Let $\Theta$ be a symmetric theta
divisor on $J$,
and
$j : X \to J$ a morphism given by
$j(x) = \omega_X - (2g-2)x$.
Then, since $j^*(\OO_J(\Theta)) = \omega_X^{\otimes 2g(g-1)}$,
we can assign the canonical adelic structure
$\overline{\omega}^{a}_X$
to $\omega_X$.
As a corollary of the above theorem,
we have the following,
which is a generalization of \cite{UlPos}.

\begin{Corollary}[cf. Corollary~\ref{cor:positive:self:int:can}]
If the polarization $\overline{B}$ is big,
then the adelic self intersection number of $\overline{\omega}^{a}_X$
with respect to $\overline{B}$ is positive, i.e.,
$\left\langle 
\overline{\omega}^{a}_X \cdot \overline{\omega}^{a}_X
\right\rangle_{\overline{B}} > 0$.
\end{Corollary}
\renewcommand{\theTheorem}{\arabic{section}.\arabic{Theorem}}

\section{Preliminaries}
For the basic notation of Arakelov Geometry, we
follow the paper \cite{MoArht}.

Let $X$ be a projective arithmetic variety with $d = \dim X_{\QQ}$,
and $\overline{L}$ a $C^{\infty}$-hermitian $\QQ$-line bundle on $X$. 
First we review several kinds of positivity of $\overline{L}$.

\medskip
$\bullet${\bf ample}:
We say $\overline{L}$ is {\em ample} if
$L$ is ample, $c_1(\overline{L})$ is a semipositive form on $X(\CC)$, and,
for a sufficiently large $n$,
$H^0(X, L^{\otimes n})$ is generated by
$\{ s \in H^0(X, L^{\otimes n}) \mid \Vert s \Vert_{\sup} < 1 \}$.

\medskip
$\bullet${\bf nef}:
We say $\overline{L}$ is {\em nef} if
$c_1(\overline{L})$ is a semipositive form on $X(\CC)$ and,
for all one-dimensional integral closed subschemes $\Gamma$ of $X$,
$\adeg \left( \rest{\overline{L}}{\Gamma} \right) \geq 0$.

\medskip
$\bullet${\bf big}:
$\overline{L}$ is said to be {\em big} if 
$\rank_{\ZZ} H^0(X, L^{\otimes m}) = O(m^d)$, and there is a non-zero
section $s$ of $H^0(X, L^{\otimes n})$ with $\Vert s \Vert_{\sup} < 1$ for
some positive integer $n$.

\medskip
$\bullet${\bf $\pmb{\QQ}$-effective}:
We say $\overline{L}$ is $\QQ$-effective,
denote by $\overline{L} \succsim 0$, if
there are a positive integer $n$ and 
a non-zero section $s \in H^0(X, L^{\otimes n})$
with $\Vert s \Vert_{\sup} \leq 1$.
Moreover, if $U$ is a non-empty Zariski open set of $X$ with
$\zero(s) \subseteq X \setminus U$, 
then we use the notation $\overline{L} \succsim_U 0$.
Let $\overline{M}$ be another $C^{\infty}$-hermitian
$\QQ$-line bundle on $X$.
If $\overline{L} \otimes \overline{M}^{\otimes -1} \succsim 0$
(resp. $\overline{L} \otimes \overline{M}^{\otimes -1} \succsim_U 0$),
then we denote this by $\overline{L} \succsim \overline{M}$
(resp. $\overline{L} \succsim_U\overline{M}$).

\begin{Proposition}
\label{prop:nef:plus:ample}
\begin{enumerate}
\renewcommand{\labelenumi}{(\arabic{enumi})}
\item
If $\overline{L}$ is a nef $C^{\infty}$-hermitian $\QQ$-line bundle and
$\overline{A}$ is an ample $C^{\infty}$-hermitian $\QQ$-line bundle, then
$\overline{L} + \epsilon \overline{A}$ is ample
for all positive rational numbers $\epsilon$.

\item
If $\overline{L}_1, \ldots, \overline{L}_{d+1}$ are nef
$C^{\infty}$-hermitian $\QQ$-line bundles, then
\[
\adeg \left( \acherncl_1(\overline{L}_1) \cdots 
\acherncl_1(\overline{L}_{d+1})
\right) \geq 0.
\]

\item
If $\overline{L}_1, \ldots, \overline{L}_{d}$ are nef
$C^{\infty}$-hermitian $\QQ$-line bundles and $\overline{M}$ is a
$\QQ$-effective $C^{\infty}$-hermitian $\QQ$-line bundle, then
\[
\adeg \left( \acherncl_1(\overline{L}_1) \cdots 
\acherncl_1(\overline{L}_{d}) \cdot \acherncl_1(\overline{M})
\right) \geq 0.
\]

\item
Let $\overline{L}_1, \ldots, \overline{L}_{d+1}$ and
$\overline{M}_1, \ldots, \overline{M}_{d+1}$ be nef $C^{\infty}$-hermitian
line bundles on $X$. 
If $\overline{M}_i \succsim \overline{L}_i$ for every $i$, then
\[
\adeg \left( \acherncl_1(\overline{M}_1) \cdots 
\acherncl_1(\overline{M}_{d+1}) \right) \geq
\adeg \left( \acherncl_1(\overline{L}_1) \cdots 
\acherncl_1(\overline{L}_{d+1}) \right).
\]
\end{enumerate}
\end{Proposition}

\Proof
(1), (2) and (3) was proved in \cite[Proposition~2.3]{MoArht}.
(4) follows from the following equation:
\begin{multline*}
\adeg \left( \acherncl_1(\overline{M}_1) \cdots 
\acherncl_1(\overline{M}_{d+1}) \right) -
\adeg \left( \acherncl_1(\overline{L}_1) \cdots 
\acherncl_1(\overline{L}_{d+1}) \right) \\
=
\sum_{i=1}^{d+1} \adeg \left( 
\acherncl_1(\overline{L}_1) \cdots \acherncl_1(\overline{L}_{i-1}) \cdot
(\acherncl_1(\overline{M}_i) - \acherncl_1(\overline{L}_i)) \cdot 
\acherncl_1(\overline{M}_{i+1}) \cdots \acherncl_1(\overline{M}_{d+1})
\right).
\end{multline*}
\QED

Moreover, the following lemmas will be used in the later sections.

\begin{Lemma}
\label{lem:exist:small:sec:zero}
Let $X$ be a projective arithmetic variety, and $\overline{L}$ a big 
$C^{\infty}$-hermitian $\QQ$-line bundle on $X$.
Let $x$ be a \rom{(}not necessarily closed\rom{)} point of $X$.
Then, there are a positive number $n$ and a non-zero section
$s \in H^0(X, L^{\otimes n})$ with $s(x) = 0$ and $\Vert s \Vert_{\sup} < 1$.
\end{Lemma}

\Proof
Since $\rank_{\ZZ} H^0(X, L^{\otimes m}) = O(m^d)$,
there are a positive number $n_0$ and 
a non-zero section
$s_0 \in H^0(X, L^{\otimes n_0})$ with $s_0(x) = 0$.
On the other hand, there is a non-zero section 
$s_1 \in H^0(X, L^{\otimes n_1})$
with $\Vert s_1 \Vert_{\sup} < 1$ for some positive integer $n_1$.
Let $n_2$ be a positive integer with 
\[
\Vert s_0 \Vert_{\sup} \Vert s_1 \Vert_{\sup}^{n_2} < 1.
\]
Thus, if we set $s = s_0 \otimes s_1^{\otimes n_2} \in
H^0(X, L^{\otimes n_0 + n_1n_2})$, then we have the desired assertion.
\QED

\begin{Lemma}
\label{lem:exist:model:ample}
Let $B$ be a projective arithmetic variety and $K$ the function field of $B$.
Let $X$ be a projective variety over $K$, and
$L$ an ample line bundle on $X$. 
Then, there are a projective arithmetic variety
$\XX$ over $B$,
and an ample $C^{\infty}$-hermitian $\QQ$-line bundle $\overline{\LL}$ on $\XX$
such that $X$ is the generic fiber of $\XX \to B$ and
$\LL$ coincides with $L$ in $\Pic(X) \otimes \QQ$.
\end{Lemma}

\Proof
Choose a sufficiently large integer $n$ such that
$\phi_{\vert L^{\otimes n} \vert}$ gives rise to
an embedding $X \hookrightarrow \PP^N_K$. Let $\XX$ be the Zariski closure
of $X$ in $\PP^N_B = \PP^N \times B$.
Since $\OO_{\PP^N}(1)$ is relative ample, 
there is an ample line bundle $Q$ on $B$
such that $\mathcal{A} = \OO_{\PP^N}(1) \otimes \pi^*(Q)$ is ample, where
$\pi$ is the natural projection $\PP^N_B \to B$.
We choose a $C^{\infty}$-hermitian metric of $\mathcal{A}$ such that
$\overline{\mathcal{A}} = (\mathcal{A}, \Vert\cdot\Vert)$ is ample.
Thus, if we set $\overline{\LL} = 
\left(\rest{\overline{\mathcal{A}}}{\XX}\right)^{\otimes 1/n}$,
then we have our assertion.
\QED

Next, let us consider the following relative positivity.

$\bullet${\bf $\pmb{\pi}$-nef} (nef with respect to a morphism):
Let $\pi : X \to B$ be a morphism of projective arithmetic varieties, and
$\overline{L}$ a $C^{\infty}$-hermitian $\QQ$-line bundle on $X$.
We say $\overline{L}$ is {\em nef with respect to $X \to B$}
(or {\em $\pi$-nef}) if the following properties
are satisfied:
\begin{enumerate}
\renewcommand{\labelenumi}{(\arabic{enumi})}
\item
For any analytic maps $h : M \to X(\CC)$ from a complex manifold $M$
to $X(\CC)$ with $\pi(h(M))$ being a point,
$c_1(h^*(\overline{L}))$ is semipositive.

\item
For every $b \in B$, the restriction $\rest{L}{X_{\overline{b}}}$ of $L$
to the geometric fiber over $b$ is nef.
\end{enumerate}
Then, we have the following lemma.

\begin{Lemma}
\label{lem:non:negative:pi:nef}
Let $\pi : X \to B$ be a morphism of projective arithmetic varieties with
$d = \dim B_{\QQ}$ and $e = \dim(X/B)$.
Let $\overline{H}_1, \ldots, \overline{H}_d$ be nef
$C^{\infty}$-hermitian $\QQ$-line bundles on $B$.
Then, we have the following.
\begin{enumerate}
\renewcommand{\labelenumi}{(\arabic{enumi})}
\item
Let $\overline{L}_1, \ldots, \overline{L}_{e}$ be $\pi$-nef
$C^{\infty}$-hermitian $\QQ$-line bundles
on $X$, and
$\overline{L}$ a $C^{\infty}$-hermitian $\QQ$-line bundle on $X$.
If there is a non-empty Zariski open set $U$ of $B$ with
$\overline{L} \succsim_{\pi^{-1}(U)} 0$, then
\[
\adeg\left(\acherncl_1(\overline{L}_1) \cdots \acherncl_1(\overline{L}_e)
\cdot \acherncl_1(\overline{L}) \cdot \acherncl_1(\pi^* \overline{H}_1)
\cdots \acherncl_1(\pi^* \overline{H}_d) \right) \geq 0.
\]

\item
Let $\overline{L}_1, \ldots, \overline{L}_{e+1}$ and
$\overline{L}'_1, \ldots, \overline{L}'_{e+1}$ be $\pi$-nef
$C^{\infty}$-hermitian $\QQ$-line bundles on $X$.
If there is a Zariski open set $U$ of $B$ such that
$\overline{L}_i \succsim_{\pi^{-1}(U)} \overline{L}'_i$ for all $i$,
then
\begin{multline*}
\adeg\left(\acherncl_1(\overline{L}_1) \cdots \acherncl_1(\overline{L}_{e+1}) 
\cdot \acherncl_1(\pi^*\overline{H}_1)
\cdots \acherncl_1(\pi^* \overline{H}_d) \right) \\
\geq
\adeg\left(\acherncl_1(\overline{L}'_1) \cdots \acherncl_1(\overline{L}'_{e+1}) 
\cdot \acherncl_1(\pi^*\overline{H}_1)
\cdots \acherncl_1(\pi^* \overline{H}_d) \right).
\end{multline*}
\end{enumerate}
\end{Lemma}

\Proof
(1)
By our assumption, there are a positive integer $n$ and a non-zero section
$s \in H^0(X, L^{\otimes n})$ such that
$\Vert s \Vert_{\sup} \leq 1$ and $\Supp(\zero(s)) \subseteq X \setminus \pi^{-1}(U)$.
Let $\zero(s) = a_1 \Delta_1 + \cdots + a_r \Delta_r$ be
the decomposition as cycles.
Then,
\addtocounter{Claim}{1}
\begin{multline}
\label{lem:non:negative:pi:nef:eqn}
n \adeg\left(\acherncl_1(\overline{L}_1) \cdots \acherncl_1(\overline{L}_e)
\cdot \acherncl_1(\overline{L}) \cdot \acherncl_1(\pi^* \overline{H}_1)
\cdots \acherncl_1(\pi^* \overline{H}_d) \right) \\
= \sum_{i=1}^r a_i \adeg\left(\acherncl_1(\rest{\overline{L}_1}{\Delta_i}) \cdots
\acherncl_1(\rest{\overline{L}_e}{\Delta_i})
\cdot \acherncl_1(\rest{\pi^* \overline{H}_1}{\Delta_i})
\cdots \acherncl_1(\rest{\pi^* \overline{H}_d)}{\Delta_i} \right) \\
+ \int_{X(\CC)} -\log(\Vert s \Vert) c_1(\overline{L}_1) \wedge \cdots \wedge
c_1(\overline{L}_e) \wedge c_1(\pi^*\overline{H}_1) \wedge \cdots \wedge
c_1(\pi^*\overline{H}_d).
\end{multline}
First, by the Fubini's theorem,
\begin{multline*}
\int_{X(\CC)} -\log(\Vert s \Vert) c_1(\overline{L}_1) \wedge \cdots \wedge
c_1(\overline{L}_e) \wedge c_1(\pi^*\overline{H}_1) \wedge \cdots \wedge
c_1(\pi^*\overline{H}_d) \\
= \int_{B(\CC)} \left( \int_{X(\CC)/B(\CC)} -\log(\Vert s \Vert) c_1(\overline{L}_1)
\wedge \cdots \wedge
c_1(\overline{L}_e) \right)  c_1(\overline{H}_1) \wedge \cdots \wedge
c_1(\overline{H}_d).
\end{multline*}
Here, by the property (1) of ``$\pi$-nef'',
\[
\int_{X(\CC)/B(\CC)} -\log(\Vert s \Vert) c_1(\overline{L}_1)
\wedge \cdots \wedge
c_1(\overline{L}_e)
\]
is a non-negative locally integrable function on $B(\CC)$.
Thus, the integral part of \eqref{lem:non:negative:pi:nef:eqn} is non-negative.
Let $b_i$ be the generic point of $\pi(\Delta_i)$.
Then, by the projection formula, we can see
\begin{multline*}
\adeg\left(\acherncl_1(\rest{\overline{L}_1}{\Delta_i}) \cdots
\acherncl_1(\rest{\overline{L}_e}{\Delta_i})
\cdot \acherncl_1(\rest{\pi^* \overline{H}_1}{\Delta_i})
\cdots \acherncl_1(\rest{\pi^* \overline{H}_d}{\Delta_i}) \right) \\
= \begin{cases}
0 & \text{if $\codim(\pi(\Delta_i)) \geq 2$} \\
\deg (\rest{L_1}{(\Delta_i)_{\bar{b}_i}} \cdots \rest{L_e}{(\Delta_i)_{\bar{b}_i}})
\adeg \left(
\acherncl_1(\rest{\overline{H}_1}{\pi(\Delta_i)})
\cdots \acherncl_1(\rest{\overline{H}_d}{\pi(\Delta_i)}) \right)
& \text{if $\codim(\pi(\Delta_i)) = 1$}
\end{cases}
\end{multline*}
Therefore, we get (1) because
\[
\deg (\rest{L_1}{(\Delta_i)_{\bar{b}_i}} \cdots \rest{L_e}{(\Delta_i)_{\bar{b}_i}}) \geq 0
\quad\text{and}\quad
\adeg \left(
\acherncl_1(\rest{\overline{H}_1}{\pi(\Delta_i)})
\cdots \acherncl_1(\rest{\overline{H}_d}{\pi(\Delta_i)}) \right) \geq 0.
\].

\medskip
(2) Since
\begin{multline*}
\acherncl_1(\overline{L}_1) \cdots \acherncl_1(\overline{L}_{e+1}) -
\acherncl_1(\overline{L}'_1) \cdots \acherncl_1(\overline{L}'_{e+1}) \\
= \sum_{i=1}^{e+1}
\acherncl_1(\overline{L}_1) \cdots \acherncl_1(\overline{L}_{i-1})
\cdot \left( \acherncl_1(\overline{L}_i) - \acherncl_1(\overline{L}'_i) \right)
\cdot \acherncl_1(\overline{L}'_{i+1}) \cdots \acherncl_1(\overline{L}'_{e+1}),
\end{multline*}
(2) is a consequence of (1).
\QED

Finally, let us consider the following lemma.

\begin{Lemma}
\label{lem:estimate:locally:int:func:below}
Let $\pi : X \to B$ be a morphism of projective arithmetic varieties, and
$\overline{L}$ a $C^{\infty}$-hermitian
line bundle on $X$. Let $U$ be a non-empty Zariski open set of $B$
such that $B \setminus U = \Supp(D)$ for some effective Cartier divisor
$D$ on $B$.
If there is a non-zero rational section $s$ of $L$ with
$\Supp(\zero(s)) \subseteq X \setminus \pi^{-1}(U)$, 
then there are a positive integer $n$ and a 
$C^{\infty}$-metric $\Vert\cdot\Vert_{nD}$
of $\OO_B(nD)$ with
\[
\pi^*(\OO_B(nD), \Vert\cdot\Vert_{nD})^{\otimes -1} \precsim_{\pi^{-1}(U)} 
\overline{L} \precsim_{\pi^{-1}(U)} \pi^*(\OO_B(nD), \Vert\cdot\Vert_{nD}).
\]
Moreover, if $D$ is ample, then we can choose $\Vert\cdot\Vert_{nD}$ such that
$(\OO_B(nD), \Vert\cdot\Vert_{nD})$ is ample.
\end{Lemma}

\Proof
First, we fix a hermitian metric $\Vert\cdot\Vert_D$ of $\OO_B(D)$.
If $D$ is ample, then we choose $\Vert\cdot\Vert_D$ such that
$(\OO_B(D), \Vert\cdot\Vert_D)$ is ample.
Find a positive integer $n$ with
\[
 -n f^*(D) \leq \zero(s) \leq n f^*(D).
\]
Let $l$ be a section of $O_Y(nD)$ with $\zero(l) = nD$.
We set $t_1 = l \otimes s^{-1}$ and $t_2 = l \otimes s$.
Then, $t_1$ and $t_2$ are global sections of
$\OO_X(nf^*(D)) \otimes L^{-1}$ and $\OO_X(nf^*(D)) \otimes L$
respectively.
Choose a sufficiently small positive number $c$ such that
if we give a norm of $\OO_B(nD)$ by $c\Vert\cdot\Vert_D^n$, then
$\Vert t_1 \Vert_{\sup} \leq 1$ and $\Vert t_2 \Vert_{\sup} \leq 1$.
Thus we get our lemma.
\QED

\section{Arithmetic height of subvarieties}
Let $K$ be a finitely generated field over $\QQ$
with $d = \trdeg_{\QQ}(K)$, 
and $\overline{B} = (B; \overline{H}_1, \ldots, \overline{H}_d)$
a polarization of $K$.
Let $X$ be a projective variety over $K$, and
$L$ a nef line bundle on $X$.
Let $\XX$ be a projective arithmetic variety over $B$ such that
$X$ is the generic fiber of $\XX \to B$, and
let $\overline{\LL}$ be a $C^{\infty}$-hermitian $\QQ$-line bundle on $\XX$
such that $\LL$ coincides with $L$ in $\Pic(X) \otimes \QQ$.
The pair $(\XX, \overline{\LL})$ is called a $C^{\infty}$-model of 
$(X, L)$. We assume that $\overline{\LL}$ is nef with respect to $\XX \to B$.
Note that if $L$ is ample, then
there is a $C^{\infty}$-model
$(\XX, \overline{\LL})$ of $(X, L)$ such that
$\overline{\LL}$ is ample by Lemma~\ref{lem:exist:model:ample}.

Let $Y$ be a subvariety of $X_{\overline{K}}$.
We assume that $Y$ is defined over a finite extension field $K'$ of $K$.
Let $B^{K'}$ be
the normalization of $B$ in $K'$, and
let $\rho^{K'} : B^{K'} \to B$ be the induced morphism.
Let $\XX^{K'}$ be the main component of $\XX \times_{B} B^{K'}$. 
We set the induced morphisms
as follows.
\[
\begin{CD}
\XX @<{\tau^{K'}}<< \XX^{K'} \\
@V{\pi}VV @VV{\pi^{K'}}V \\
B     @<{\rho^{K'}}<< B^{K'}
\end{CD}
\]
Let $\YY$ be the Zariski closure of $Y$ in $\XX^{K'}$.
Then the naive height $h^{\overline{B}}_{(\XX, \overline{\LL})}(Y)$ of $Y$
with respect to $(\XX, \overline{\LL})$ and $\overline{B}$
is defined by
\begin{multline*}
h^{\overline{B}}_{(\XX, \overline{\LL})}(Y) \\
= \frac{\adeg\left( 
\acherncl_1\left(\rest{{\tau^{K'}}^*(\overline{\LL})}{\YY}\right)^{\cdot \dim Y  +1} \cdot 
\acherncl_1\left(\rest{{\pi^{K'}}^*({\rho^{K'}}^*(\overline{H}_1))}{\YY}\right) \cdots
\acherncl_1\left(\rest{{\pi^{K'}}^*({\rho^{K'}}^*(\overline{H}_d))}{\YY}\right)
\right)}{[K':K](\dim Y + 1)\deg(\rest{L}{Y}^{\dim Y})}.
\end{multline*}
Note that the above definition does not depend on the choice of $K'$
by the projection formula.
Here we have the following proposition.
By this proposition, we may denote by $h^{\overline{B}}_L$
the class of $h^{\overline{B}}_{(\XX, \overline{\LL})}$ modulo
the set of bounded functions.
Moreover, we say $h^{\overline{B}}_L$ is the height function
associated with $L$ and $\overline{B}$.

\begin{Proposition}
\label{prop:bound:diff:two:model}
Let $(\XX', \overline{\LL}')$ be another model of $(X, L)$ over $B$
such that $\overline{\LL}'$ is nef with respect to $\XX' \to B$.
Then, there is a constant $C$ such that
\[
\left\vert 
h^{\overline{B}}_{(\XX, \overline{\LL})}(Y) - h^{\overline{B}}_{(\XX', \overline{\LL}')}(Y)
\right\vert \leq C
\]
for all subvarieties $Y$ of $X_{\overline{K}}$.
\end{Proposition}

\Proof
Let $U$ be a Zariski open set of $B$ with
$\XX_U = \XX'_U$ and $\LL_U = \LL'_U$ in $\Pic(\XX_U) \otimes \QQ$.
Let $A$ be an ample line bundle on $B$ and $I$ the defining ideal of
$B \setminus U$. Then, there is a non-zero section $t$ of 
$H^0(B, A^{\otimes m} \otimes I)$
for some positive integer $m$. Thus, $B \setminus U \subseteq \Supp(\zero(s))$.
Therefore, shrinking $U$, we may assume that 
there is an effective ample Cartier divisor
$D$ on $B$ with $\Supp(D) =  B \setminus U$.

Let $\mu : \ZZZ \to \XX$ and $\mu' : \ZZZ \to \XX'$ be birational morphisms of
projective arithmetic varieties such that $\mu$ and $\mu'$ are the identity map over $\XX_U$. 
Then, $h_{(\XX, \overline{\LL})}^{\overline{B}}(Y) = 
h_{(\ZZZ, \mu^*(\overline{\LL}))}^{\overline{B}}(Y)$ and
$h_{(\XX', \overline{\LL}')}^{\overline{B}}(Y) = 
h_{(\ZZZ, {\mu'}^*(\overline{\LL}'))}^{\overline{B}}(Y)$
for all subvarieties $Y$ of $X_{\overline{K}}$.
Thus, to prove our proposition,
we may assume that $\XX = \XX'$.

First of all, by Lemma~\ref{lem:estimate:locally:int:func:below},
there is a nef $C^{\infty}$-hermitian line bundle $\overline{T}$ on $B$ such that
\addtocounter{Claim}{1}
\begin{equation}
\label{prop:bound:diff:two:model:eqn:1}
\pi^*(\overline{T})^{\otimes -1} \precsim_{\pi^{-1}(U)} 
\overline{\LL} \otimes {\overline{\LL}'}^{\otimes -1}
\precsim_{\pi^{-1}(U)} \pi^*(\overline{T}),
\end{equation}
where $\pi : \XX \to B$ is the canonical morphism.
Let $Y$ be a subvariety of $X_{\overline{K}}$.
We assume that $Y$ is defined over a finite extension field $K'$ of $K$.
Let $B^{K'}$ be
the normalization of $B$ in $K'$, and
$\XX^{K'}$ the main components of $\XX \times_{B} B^{K'}$. 
Let $\YY$ be the closure of $Y$ in $\XX^{K'}$. Then,
\[
h^{\overline{B}}_{(\XX, \overline{\LL})}(Y)
= \frac{\adeg\left( 
\acherncl_1\left(\rest{\overline{\LL}^{K'}}{\YY}\right)^{\cdot \dim Y  +1} \cdot 
\acherncl_1\left(\rest{\overline{H}_1^{K'}}{\YY}\right) \cdots
\acherncl_1\left(\rest{\overline{H}_d^{K'}}{\YY}\right)
\right)}{[K':K](\dim Y + 1)\deg(\rest{L}{Y}^{\dim Y})}
\]
and
\[
h^{\overline{B}}_{(\XX', \overline{\LL}')}(Y)
= \frac{\adeg\left( 
\acherncl_1\left(\rest{{\overline{\LL}'}^{K'}}{\YY}\right)^{\cdot \dim Y  +1} \cdot 
\acherncl_1\left(\rest{\overline{H}_1^{K'}}{\YY}\right) \cdots
\acherncl_1\left(\rest{\overline{H}_d^{K'}}{\YY}\right)
\right)}{[K':K](\dim Y + 1)\deg(\rest{L}{Y}^{\dim Y})},
\]
where $\overline{\LL}^{K'}$,
${\overline{\LL}'}^{K'}$ and
$\overline{H}_i^{K'}$'s are
pullbacks of
$\overline{\LL}$,
$\overline{\LL}'$ and
$\overline{H}_i$'s
to $\XX^{K'}$ respectively.
Here, by virtue of \eqref{prop:bound:diff:two:model:eqn:1},
\[
{\overline{\LL}'}^{K'} \otimes {\overline{T}^{K'}}^{\otimes -1}
\precsim \overline{\LL}^{K'}
\precsim
{\overline{\LL}'}^{K'} \otimes \overline{T}^{K'}.
\]
Therefore, by (2) of Lemma~\ref{lem:non:negative:pi:nef},
we can see that
\begin{multline*}
\left\vert
\adeg\left( 
\acherncl_1\left(\rest{\overline{\LL}^{K'}}{\YY}\right)^{\cdot \dim Y  +1} \cdot 
\acherncl_1\left(\rest{\overline{H}_1^{K'}}{\YY}\right) \cdots
\acherncl_1\left(\rest{\overline{H}_d^{K'}}{\YY}\right)
\right) \right. \\
\left. -
\adeg\left( 
\acherncl_1\left(\rest{{\overline{\LL}'}^{K'}}{\YY}\right)^{\cdot \dim Y  +1} \cdot 
\acherncl_1\left(\rest{\overline{H}_1^{K'}}{\YY}\right) \cdots
\acherncl_1\left(\rest{\overline{H}_d^{K'}}{\YY}\right)
\right)
\right\vert \\
\leq
[K':K](\dim Y + 1)\deg(\rest{L}{Y}^{\dim Y}) \adeg \left(
\overline{T} \cdot \overline{H}_1 \cdots \overline{H}_d \right).
\end{multline*}
Thus we get our proposition.
\QED

\renewcommand{\theTheorem}{\arabic{section}.\arabic{subsection}.\arabic{Theorem}}
\renewcommand{\theClaim}{\arabic{section}.\arabic{subsection}.\arabic{Theorem}.\arabic{Claim}}
\renewcommand{\theequation}{\arabic{section}.\arabic{subsection}.\arabic{Theorem}.\arabic{Claim}}
\section{Adelic sequence and adelic structure}
\subsection{Adelic sequence, adelic structure and adelic line bundle}
\label{subsec:adelic:sequence}
Let $K$ be a finitely generated field over $\QQ$
with $d = \trdeg_{\QQ}(K)$, 
and $\overline{B} = (B; \overline{H}_1, \ldots, \overline{H}_d)$
a polarization of $K$.
Let $X$ be a projective variety over $K$, and
$L$ a nef line bundle on $X$.

A sequence of $C^{\infty}$-models $\{ (\XX_n, \overline{\LL}_n) \}$
of $(X, L)$ is called an {\em adelic sequence of $(X, L)$} (with
respect to $\overline{B}$) if 
$\overline{\LL}_n$ is nef with respect to $\XX_n \to B$
for every $n$, and there is a non-empty
Zariski open set $U$ of $B$ with
following properties:
\begin{enumerate}
\renewcommand{\labelenumi}{(\arabic{enumi})}
\item
$\rest{\XX_n}{U} = \rest{\XX_m}{U}$ (say $\XX_U$) and 
$\rest{\LL_n}{U} = \rest{\LL_m}{U}$
in $\Pic(\XX_U) \otimes \QQ$
for all $n, m$.

\item
For each $n, m$,
there are a projective arithmetic variety $\XX_{n, m}$ over $B$,
birational morphisms $\mu_{n, m}^{n} : \XX_{n,m} \to \XX_n$ and 
$\mu_{n, m}^{m} : \XX_{n, m} \to \XX_m$, and
a nef $C^{\infty}$-hermitian $\QQ$-line bundle $\overline{D}_{n, m}$ on $B$
such that
\[
 \pi_{n, m}^*(\overline{D}_{n, m}^{\otimes -1}) \precsim_{\pi_{n, m}^{-1}(U)}
 (\mu_{n, m}^{n})^*(\overline{\LL}_n) \otimes 
 (\mu_{n, m}^{m})^*(\overline{\LL}_m^{\otimes -1}) 
\precsim_{\pi_{n, m}^{-1}(U)}
 \pi_{n, m}^*(\overline{D}_{n, m})
\]
and that
\[
\adeg\left( \acherncl_1(\overline{D}_{n, m})
\cdot \acherncl_1(\overline{H}_1) \cdots \acherncl_1(\overline{H}_d)
\right) \longrightarrow 0
\]
as $n, m \longrightarrow \infty$,
where $\pi_{n,m}$ is the natural morphism $\XX_{n,m} \to B$.
\end{enumerate}
The open set $U$ as above is called a {\em common base} of the sequence
$\{ (\XX_n, \overline{\LL}_n) \}$.
Note that if $U'$ is a non-empty Zariski open set of $U$, then
$U'$ is also a common base of $\{ (\XX_n, \overline{\LL}_n) \}$.

Let $\{ (\YY_n, \overline{\MM}_n) \}$ be another adelic sequence
of $(X, L)$. We say $\{ (\XX_n, \overline{\LL}_n) \}$ 
is {\em equivalent} to $\{ (\YY_n, \overline{\MM}_n) \}$,
denoted by $\{ (\XX_n, \overline{\LL}_n) \} \sim
\{ (\YY_n, \overline{\MM}_n) \}$,
if the concatenated sequence
\[
 (\XX_1, \overline{\LL}_1), (\YY_1, \overline{\MM}_1),
\cdots, (\XX_n, \overline{\LL}_n), (\YY_n, \overline{\MM}_n), \cdots
\]
is adelic.  In other words, if we choose a suitable common base
$U$, then, for each $n$, there are a projective arithmetic variety $\ZZZ_{n}$
over $B$, birational morphisms $\mu_{n} : \ZZZ_{n} \to \XX_n$ and 
$\nu_{n} : \ZZZ_{n} \to \YY_n$, and
a nef $C^{\infty}$-hermitian $\QQ$-line bundle $\overline{D}_{n}$ on $B$
such that
\[
 \pi_{\ZZZ_n}^*(\overline{D}_{n}^{\otimes -1}) 
 \precsim_{\pi^{-1}_{\ZZZ_n}(U)}
 \mu_{n}^*(\overline{\LL}_n) \otimes 
 \nu_{n}^*(\overline{\MM}_n^{\otimes -1}) 
 \precsim_{\pi^{-1}_{\ZZZ_n}(U)}
 \pi_{\ZZZ_n}^*(\overline{D}_{n})
\]
and that
\[
\lim_{n \to \infty} \adeg\left( \acherncl_1(\overline{D}_{n})
\cdot \acherncl_1(\overline{H}_1) \cdots \acherncl_1(\overline{H}_d)
\right) = 0,
\]
where $\pi_{\ZZZ_n}$ is the natural morphism $\ZZZ_{n} \to B$.

An equivalent class of adelic sequences of $(X, L)$ is called
an {\em adelic structure of $L$} (with respect to $\overline{B}$).
Further, a line bundle $L$ with an adelic structure is called
an {\em adelic line bundle} and is often denoted by $\overline{L}$
for simplicity. 
If an adelic line bundle $\overline{L}$ is given by
an adelic sequence $\{ (\XX_n, \overline{\LL}_n) \}$,
then we denote this by $\overline{L} = \lim_{n \to \infty}
(\XX_n, \overline{\LL}_n)$.
Moreover, we say $\overline{L}$ is nef
if $\overline{L} = \lim_{n\to\infty}
(\XX_n, \overline{\LL}_n)$ and
$\overline{\LL}_n$ is nef for $n \gg 0$.

Let $g : Y \to X$ be a morphism
of projective varieties over $K$, and
$\overline{L}$ an adelic line bundle on $X$.
We assume that $\overline{L}$ is given by
an adelic sequence $\{ (\XX_n, \overline{\LL}_n) \}$.
Let us fix a morphism $g_n : \YY_n \to \XX_n$ of
projective arithmetic varieties over $B$ for each $n$
with the following properties:
\begin{enumerate}
\renewcommand{\labelenumi}{(\alph{enumi})}
\item
$g_n : \YY_n \to \XX_n$ coincides with $g : Y \to X$ over $K$ for every $n$.

\item
There is a non-empty Zariski open set $U$ of $B$
such that $\rest{\YY_n}{U} = \rest{\YY_m}{U}$,
$\rest{\XX_n}{U} = \rest{\XX_m}{U}$, and
$\rest{g_n}{U} = \rest{g_m}{U}$ for all $n, m$
\end{enumerate}
Then it is not difficult to see that
$\{ (\YY_n, g_n^*(\overline{\LL}_n)) \}$
is an adelic sequence of $(Y, g^*(L))$.
We denote by $g^*(\overline{L})$ the adelic structure given by
$\{ (\YY_n, g_n^*(\overline{\LL}_n)) \}$.
Note that this adelic structure does not depend on the choice
of the adelic sequence $\{ (\XX_n, \overline{\LL}_n) \}$ and
the morphisms $g_n : \YY_n \to \XX_n$.

\setcounter{Theorem}{0}
\subsection{Adelic sequence by an endomorphism}
Let $K$ be a finitely generated field over $\QQ$
with $d = \trdeg_{\QQ}(K)$, 
and $\overline{B} = (B; \overline{H}_1, \ldots, \overline{H}_d)$
a polarization of $K$.
Let $X$ be a projective variety over $K$, and $L$ an ample line bundle on $X$.
We assume that there is a surjective morphism $f : X \to X$
and an integer $d \geq 2$
with $L^{\otimes d} \simeq f^*(L)$.
Let $(\XX, \overline{\LL})$ be a $C^{\infty}$-model of $(X, L)$ such that
$\overline{\LL}$ is nef with respect to $\XX \to B$.
Note that the existence of a $C^{\infty}$-model
$(\XX, \overline{\LL})$ of $(X, L)$ with $\LL$ being nef
with respect to $\XX \to B$
is guaranteed by Lemma~\ref{lem:exist:model:ample}.
Then, there is a Zariski open set $U$ of $B$ such that
$f$ extends to $f_U: \XX_U \to \XX_U$ and
$\LL_U^{\otimes d} = f_U^*(\LL_U)$ in $\Pic(\XX_U) \otimes \QQ$.
Let $\XX_n$ be the normalization of 
$\XX_U \overset{f_U^n}{\longrightarrow} \XX_U \to \XX$, and
$f_n : \XX_n \to \XX$ the induced morphism.
Then, we have the following proposition.

\begin{Proposition}
\label{prop:f:adelic:sequence}
\begin{enumerate}
\renewcommand{\labelenumi}{(\arabic{enumi})}
\item
$\left\{ (\XX_n, f_n^*(\overline{\LL})^{\otimes d^{-n}}) \right\}$
is an adelic sequence
of $(X, L)$.
Moreover, if $\overline{\LL}$ is nef, then
the adelic line bundle
$\lim_{n\to\infty}(\XX_n, f_n^*(\overline{\LL})^{\otimes d^{-n}})$
is nef.

\item
Let $f' : X \to X$ be another surjective morphism
with 
$L^{\otimes d'} \simeq {f'}^*(L)$
for some $d' \geq 2$. 
Let $(\XX', \overline{\LL}')$ be another 
$C^{\infty}$-model of $(X, L)$ such that
$\overline{\LL}'$ is nef with respect to $\XX' \to B$.
Let $U'$ be a non-empty Zariski open set of $B$ such that
$f'$ extends to $f'_{U'}: \XX'_{U'} \to \XX'_{U'}$ and
${\LL'_{U'}}^{\otimes d'} ={f'_{U'}}^*(\LL'_{U'})$
in $\Pic(\XX'_{U'}) \otimes \QQ$.
Let $\XX'_n$ be the normalization of 
$\XX'_{U'} \overset{{f'_{U'}}^n}{\longrightarrow} \XX'_{U'} \to \XX'$, and
$f'_n : \XX'_n \to \XX'$ the induced morphism.
If $f \cdot f' = f' \cdot f$, then
\[
\left\{ (\XX_n, f_n^*(\overline{\LL})^{\otimes d^{-n}}) \right\} \sim
\left\{ (\XX'_n, {f'_n}^*(\overline{\LL}')^{\otimes {d'}^{-n}}) \right\}.
\]
\end{enumerate}
\end{Proposition}

\begin{Definition}[$f$-adelic structure]
The adelic sequence
$\left\{ (\XX_n, f_n^*(\overline{\LL})^{\otimes d^{-n}}) \right\}$
in the above proposition gives rise to
the adelic structure on $L$, which
is called the {\em $f$-adelic structure of $L$}.
The line bundle $L$
with this  adelic structure is denoted by $\overline{L}^f$,
i.e.,
$\overline{L}^f = \lim_{n\to\infty} 
(\XX_n, f_n^*(\overline{\LL})^{\otimes d^{-n}})$.
Considering a case ``$f = f'$'' in (2),
we can see that $\overline{L}^f$ does not depend on the choice
of the $C^{\infty}$-model $(\XX, \overline{\LL})$.
Moreover, (2) says us that
if $f \cdot f' = f' \cdot f$, then
$\overline{L}^f = \overline{L}^{f'}$.
Further, $\overline{L}^f$ is nef by the second assertion
of (1) and Lemma~\ref{lem:exist:model:ample}.
\end{Definition}

{\sl Proof of Proposition}~\ref{prop:f:adelic:sequence}.\quad
In the same way as in the proof of
Proposition~\ref{prop:bound:diff:two:model}, 
shrinking $U$ if necessarily, we may assume that 
there is an effective ample Cartier divisor
$D$ on $B$ with $\Supp(D) =  B \setminus U$.

\medskip
(1) For simplicity, we denote $f_n^*(\overline{\LL})^{\otimes d^{-n}}$ 
by $\overline{\LL}_n$.
From now on, we treat the group structure of the Picard group additively.
Note that $\XX_0 = \XX$ and $\overline{\LL}_0 = \overline{\LL}$.
Let $\YY$ be a projective arithmetic variety over $B$ such that
there are birational morphisms 
$\rho_0 : \YY \to \XX_0$ and $\rho_1 : \YY \to \XX_1$,
which are the identity map over $U$.
We fix $n > m \geq 0$.
Let $\ZZZ$ be a projective arithmetic variety over $B$ with the following
properties:
\begin{enumerate}
\renewcommand{\labelenumi}{(\alph{enumi})}
\item
$\ZZZ_U = \XX_U$.

\item
For each $m \leq i \leq n$, there is a birational morphism
$\mu_i : \ZZZ \to \XX_i$, which is the identity map over $U$.

\item
For each $m \leq j < n$, there is a morphism
$g_j : \ZZZ \to \YY$ which is an extension of $f_U^j : \ZZZ_U \to \YY_U$.
\end{enumerate}
Here we claim the following.
\renewcommand{\theClaim}{\arabic{section}.\arabic{subsection}.\arabic{Theorem}}
\addtocounter{Theorem}{1}
\begin{Claim}
\label{claim:prop:f:adelic:sequence:1}
\begin{enumerate}
\renewcommand{\labelenumi}{(\roman{enumi})}
\item
$\mu_{j+1}^*(\overline{\LL}_{j+1}) = 
d^{-j} g_{j}^*(\rho_1^*(\overline{\LL}_1))$
for each $m \leq j < n$.

\item
$\mu_j^*(\overline{\LL}_j) = 
d^{-j} g_{j}^*(\rho_0^*(\overline{\LL}_0))$
for each $m \leq j < n$.

\item
${\displaystyle \mu_n^*(\overline{\LL}_n) - \mu_m^*(\overline{\LL}_m) =
\sum_{j=m}^{n-1} d^{-j} g_j^*(\rho_1^*(\overline{\LL}_1) - 
\rho_0^*(\overline{\LL}_0))}$.
\end{enumerate}
\end{Claim}
\renewcommand{\theClaim}{\arabic{section}.\arabic{subsection}.\arabic{Theorem}.\arabic{Claim}}

(i)
Let us consider the following two morphisms between $\ZZZ$ and $\XX_0$:
\[
\ZZZ \overset{g_{j}}{\longrightarrow} \YY \overset{\rho_1}{\longrightarrow}
\XX_1 \overset{f_1}{\longrightarrow} \XX_0
\qquad\text{and}\qquad
\ZZZ \overset{\mu_{j+1}}{\longrightarrow} 
\XX_{j+1} \overset{f_{j+1}}{\longrightarrow} \XX_0.
\]
These are same over $U$. Thus, so are over $B$.
Therefore, $g_{j}^* \rho_1^* f_1^*(\overline{\LL}) =
\mu_{j+1}^* f_{j+1}^*(\overline{\LL})$, which shows us the assertion of (i).

\smallskip
(ii) In the same way as above, we can see
$\mu_j \cdot f_j = \rho_0 \cdot g_j$. Thus we get (ii).

\smallskip
(iii)
Since ${\displaystyle
\mu_n^*(\overline{\LL}_n) - \mu_m^*(\overline{\LL}_m) =
\sum_{j=m}^{n-1} \mu_{j+1}^*(\overline{\LL}_{j+1}) - 
\mu_j^*(\overline{\LL}_j),}$
this is a consequence of (i) and (ii).

\bigskip
By Lemma~\ref{lem:estimate:locally:int:func:below},
there is an ample $C^{\infty}$-hermitian line bundle
$\overline{\Delta}$ on $B$ such that
\[
-\pi_{\YY}^*(\overline{\Delta}) \precsim_{\pi_{\YY}^{-1}(U)}
\rho_1^*(\overline{\LL}_1) - \rho_0^*(\overline{\LL}_0) 
\precsim_{\pi_{\YY}^{-1}(U)}
\pi_{\YY}^*(\overline{\Delta}).
\]
Hence, by (iii) of the above claim, we get
\[
- \left( \sum_{j=m}^{n-1} d^{-j} \right) \pi_{\ZZZ}^*(\overline{\Delta})
\precsim_{\pi_{\ZZZ}^{-1}(U)}
\mu_n^*(\overline{\LL}_n) - \mu_m^*(\overline{\LL}_m) \precsim_{\pi_{\ZZZ}^{-1}(U)}
\left( \sum_{j=m}^{n-1} d^{-j} \right) \pi_{\ZZZ}^*(\overline{\Delta}).
\]
Thus, we obtain the first assertion of (1).
The second assertion is obvious.

\bigskip
(2) Let us consider the following cases:
\par\smallskip
Case 1 : $f = f'$.
\par
Case 2 : $\XX = \XX'$ and $\overline{\LL} = \overline{\LL}'$.
\par\noindent\smallskip
Clearly, it is sufficient to check (2) under the assumption Case 1 or Case 2.

\medskip
{\bf Case 1} : In this case, we assume $f = f'$.
Shrinking $U$ and $U'$, we may assume that $U = U'$, $\XX_U = \XX'_{U'}$ and
$\LL_U = \LL'_{U'}$ in $\Pic(\XX_U) \otimes \QQ$.
For each $n \geq 0$,
let $\ZZZ_n$ be a projective arithmetic variety over $B$ such that
there are birational morphisms 
$\nu_n : \ZZZ_n \to \XX_n$ and $\nu'_n : \ZZZ_n \to \XX'_n$, which are
the identity map over $U$.
We may assume that there is a morphism
$g_n : \ZZZ_n \to \ZZZ_0$ such that the following diagrams are commutative:
\[
\begin{CD}
\ZZZ_0 @<{g_n}<< \ZZZ_n \\
@V{\nu_0}VV      @VV{\nu_n}V \\
\XX_0 @<{f_n}<<  \XX_n
\end{CD}
\qquad\qquad
\begin{CD}
\ZZZ_0 @<{g_n}<< \ZZZ_n \\
@V{\nu'_0}VV      @VV{\nu'_n}V \\
\XX'_0 @<{f'_n}<<  \XX'_n
\end{CD}
\]
Then,
\[
d^{-n}\nu_n^*(f_n^*(\overline{\LL})) -
d^{-n} {\nu'_n}^*({f'_n}^*(\overline{\LL}') )
= d^{-n} g_n^* \left(
\nu_0^*(\overline{\LL}) - {\nu'_0}^*(\overline{\LL}')
\right).
\]
By Lemma~\ref{lem:estimate:locally:int:func:below},
there is an ample $C^{\infty}$-hermitian line bundle
$\overline{\Delta}$ on $B$ such that
\[
-\pi_{\ZZZ_0}^*(\overline{\Delta}) \precsim_{\pi^{-1}_{\ZZZ_0}(U)}
\nu_0^*(\overline{\LL}) - {\nu'_0}^*(\overline{\LL}') 
\precsim_{\pi^{-1}_{\ZZZ_0}(U)}
\pi_{\ZZZ_0}^*(\overline{\Delta}).
\]
Therefore, we have
\[
-d^{-n} \pi_{\ZZZ_n}^*(\overline{\Delta}) \precsim_{\pi^{-1}_{\ZZZ_n}(U)}
d^{-n} \nu_n^*(f_n^*(\overline{\LL}) ) -
d^{-n} {\nu'_n}^*({f'_n}^*(\overline{\LL}'))
\precsim_{\pi^{-1}_{\ZZZ_n}(U)}
d^{-n} \pi_{\ZZZ_n}^*(\overline{\Delta}),
\]
which shows us our assertion in this case.

\medskip
{\bf Case 2} : In this case, we assume that
$\XX = \XX'$ and $\overline{\LL} = \overline{\LL}'$.
We denote
$f_n^*(\overline{\LL})^{\otimes d^{-n}}$ and
${f'_n}^*(\overline{\LL})^{\otimes {d'}^{-n}}$
by $\overline{\LL}_n$ and $\overline{\LL}'_n$
respectively.
Let $\YY$ be a projective arithmetic variety over $B$ such that
there are birational morphisms
$\rho : \YY \to \XX$,
$\rho_1 : \YY \to \XX_{1}$, and
$\rho'_1 : \YY \to \XX'_{1}$, which are the identity map over $U$.
We fix $n > 0$.
Let $\ZZZ$ be a projective arithmetic variety over $B$ with the following
properties:
\begin{enumerate}
\renewcommand{\labelenumi}{(\alph{enumi})}
\item
$\ZZZ_U = \XX_U$.

\item
For each $0 \leq i \leq n$, there are birational morphisms
$\mu_i : \ZZZ \to \XX_i$ and $\mu_i' : \ZZZ \to \XX'_i$, which are
the identity map over $U$.

\item
For each $0 \leq j \leq n$, there are morphisms
$g_j : \ZZZ \to \YY$ and
$g'_j : \ZZZ \to \YY$ which are extensions of $f_U^j : \ZZZ_U \to \YY_U$ and
${f'}_U^j : \ZZZ_U \to \YY_U$ respectively.
\end{enumerate}
Note that $\mu_0 = \mu'_0$.
Then, in the same way as in (1)
((iii) of Claim~\ref{claim:prop:f:adelic:sequence:1}), 
we can see
\renewcommand{\theequation}{\arabic{section}.\arabic{subsection}.\arabic{Theorem}}
\addtocounter{Theorem}{1}
\begin{equation}
\label{eqn:prop:f:adelic:sequence:1}
\mu_n^*(\overline{\LL}_n) - \mu_0^*(\overline{\LL}) =
\sum_{j=0}^{n-1} d^{-j} g_j^*(\rho_1^*(\overline{\LL}_1) - 
\rho^*(\overline{\LL}))
\end{equation}
and
\addtocounter{Theorem}{1}
\begin{equation}
\label{eqn:prop:f:adelic:sequence:2}
{\mu'_n}^*(\overline{\LL}'_n) - \mu_0^*(\overline{\LL}) =
\sum_{j=0}^{n-1} {d'}^{-j} {g'_j}^*({\rho'_1}^*(\overline{\LL}'_1) - 
\rho^*(\overline{\LL})).
\end{equation}
\renewcommand{\theequation}{\arabic{section}.\arabic{Theorem}.\arabic{Claim}}
Let $\ZZZ_n$ (resp. $\ZZZ'_n$) be the normalization of
$\ZZZ_U \overset{f_U^n}{\longrightarrow} \ZZZ_U \to \ZZZ$ (resp.
$\ZZZ_U \overset{{f'_U}^n}{\longrightarrow} \ZZZ_U \to \ZZZ$)
and let $h_n : \ZZZ_n \to \ZZZ$ (resp. $h'_n : \ZZZ'_n \to \ZZZ$) 
be the induced morphism.
Moreover, let $\mathcal{T}$ be 
a projective arithmetic variety over $B$ such that
there are birational morphisms 
$\tau : \mathcal{T} \to \ZZZ_n$, $\tau' : \mathcal{T} \to \ZZZ'_n$,
$\sigma : \mathcal{T} \to \XX_{n}$, and
$\sigma' : \mathcal{T} \to \XX'_{n}$,
which are the identity map over $U$.
Now we have a lot of morphisms, so that we summarize them.
The following morphisms are birational and the identity map over $U$.
\[
\begin{CD}
\YY @>{\rho}>> \XX
\end{CD}
\qquad
\begin{CD}
\YY @>{\rho_1}>> \XX_1 \\
\YY @>{\rho'_1}>> \XX'_1
\end{CD}
\qquad
\begin{CD}
\ZZZ @>{\mu_i}>> \XX_i \\
\ZZZ @>{\mu'_i}>> \XX'_i
\end{CD}
\qquad
\begin{CD}
\mathcal{T} @>{\tau}>> \ZZZ_n \\
\mathcal{T} @>{\tau'}>> \ZZZ'_n
\end{CD}
\qquad
\begin{CD}
\mathcal{T} @>{\sigma}>> \XX_n \\
\mathcal{T} @>{\sigma'}>> \XX'_n
\end{CD}
\]
Moreover, the following morphisms are extensions of the power of $f$ or ${f'}$.
\[
\begin{CD}
\XX_n @>{f_n}>> \XX \ \text{($f^n$ over $U$)} \\
\XX'_n @>{f'_n}>> \XX \ \text{(${f'}^n$ over $U$)}
\end{CD}
\qquad
\begin{CD}
\ZZZ @>{g_j}>> \YY \ \text{($f^j$ over $U$)} \\
\ZZZ @>{g'_j}>> \YY \ \text{(${f'}^j$ over $U$)}
\end{CD}
\qquad
\begin{CD}
\ZZZ_n @>{h_n}>> \ZZZ \ \text{($f^n$ over $U$)} \\
\ZZZ'_n @>{h'_n}>> \ZZZ \ \text{(${f'}^n$ over $U$)}
\end{CD}
\]
Here, $f_n \cdot \mu_n \cdot h'_n \cdot \tau' =
f'_n \cdot \mu'_n \cdot h_n \cdot \tau$ over $U$ because
$f \cdot f' = f' \cdot f$. Hence,
so is over $B$ as $\mathcal{T} \to \XX$. Thus,
\begin{multline*}
{d'}^{-n}{\tau'}^*{h'_n}^*\left(\mu_n^*(\overline{\LL}_n) - 
\mu_0^*(\overline{\LL})\right) -
d^{-n}\tau^*h_n^* \left( {\mu'_n}^*(\overline{\LL}'_n) - 
\mu_0^*(\overline{\LL}) \right) \\
=
d^{-n}\tau^*h_n^*\mu_0^*(\overline{\LL}) - 
{d'}^{-n}{\tau'}^*{h'_n}^*\mu_0^*(\overline{\LL})
\end{multline*}
Moreover, since  $\mu_0 \cdot h_n \cdot \tau = f_n \cdot \sigma$ and
$\mu_0 \cdot h'_n \cdot \tau' = f'_n \cdot \sigma'$, by the above equation,
we have
\renewcommand{\theequation}{\arabic{section}.\arabic{subsection}.\arabic{Theorem}}
\addtocounter{Theorem}{1}
\begin{equation}
\label{eqn:prop:f:adelic:sequence:3}
{d'}^{-n}{\tau'}^*{h'_n}^*\left(\mu_n^*(\overline{\LL}_n) - 
\mu_0^*(\overline{\LL})\right) -
d^{-n}\tau^*h_n^* \left( {\mu'_n}^*(\overline{\LL}'_n) - 
\mu_0^*(\overline{\LL}) \right) =
\sigma^*(\overline{\LL}_n) - {\sigma'}^*(\overline{\LL}'_n)
\end{equation}
\renewcommand{\theequation}{\arabic{section}.\arabic{Theorem}.\arabic{Claim}}
On the other hand, by Lemma~\ref{lem:estimate:locally:int:func:below},
we can find ample $C^{\infty}$-hermitian
$\QQ$-line bundles $\overline{\Delta}$ and
$\overline{\Delta}'$ such that
\[
-\pi_{\YY}^*(\overline{\Delta}) \precsim_{\pi^{-1}_{\YY}(U)}
\rho_1^*(\overline{\LL}_1) - 
\rho^*(\overline{\LL})
\precsim_{\pi^{-1}_{\YY}(U)} \pi_{\YY}^*(\overline{\Delta})
\]
and
\[
-\pi_{\YY}^*(\overline{\Delta}') \precsim_{\pi^{-1}_{\YY}(U)} {\rho'_1}^*(\overline{\LL}'_1) - 
\rho^*(\overline{\LL})
\precsim_{\pi^{-1}_{\YY}(U)} \pi_{\YY}^*(\overline{\Delta}').
\]
Therefore, if we set
\[
d_n = \left({d'}^{-n} \sum_{j=0}^{n-1} d^{-j} \right)
\quad\text{and}\quad
d'_n = \left(d^{-n} \sum_{j=0}^{n-1} {d'}^{-j} \right),
\]
then, by \eqref{eqn:prop:f:adelic:sequence:1} and
\eqref{eqn:prop:f:adelic:sequence:2},
\[
- d_n \pi_{\mathcal{T}}(\overline{\Delta})
\precsim_{\pi^{-1}_{\mathcal{T}}(U)}
{d'}^{-n}{\tau'}^*{h'_n}^*\left(\mu_n^*(\overline{\LL}_n) - 
\mu_0^*(\overline{\LL})\right) 
\precsim_{\pi^{-1}_{\mathcal{T}}(U)}
d_n \pi_{\mathcal{T}}(\overline{\Delta})
\]
and
\[
- d'_n
\pi_{\mathcal{T}}(\overline{\Delta}')
\precsim_{\pi^{-1}_{\mathcal{T}}(U)}
d^{-n}\tau^*h_n^* \left( {\mu'_n}^*(\overline{\LL}'_n) - 
\mu_0^*(\overline{\LL}) \right)
\precsim_{\pi^{-1}_{\mathcal{T}}(U)}
d'_n
\pi_{\mathcal{T}}(\overline{\Delta}').
\]
Hence, using \eqref{eqn:prop:f:adelic:sequence:3},
\[
- \pi_{\mathcal{T}}\left(d_n \overline{\Delta}
+ d'_n \overline{\Delta}'\right)
\precsim_{\pi^{-1}_{\mathcal{T}}(U)}
\sigma^*(\overline{\LL}_n) - {\sigma'}^*(\overline{\LL}'_n)
\precsim_{\pi^{-1}_{\mathcal{T}}(U)}
\pi_{\mathcal{T}}\left(d_n \overline{\Delta}
+ d'_n \overline{\Delta}'\right).
\]
Therefore, we have this case because
$\lim_{n\to\infty} d_n = \lim_{n\to\infty} d'_n = 0$.
\QED

\section{Adelic intersection number and adelic height}
\subsection{Adelic intersection number}
Let $K$ be a finitely generated field over $\QQ$
with $d = \trdeg_{\QQ}(K)$, 
and $\overline{B} = (B; \overline{H}_1, \ldots, \overline{H}_d)$
a polarization of $K$.

\begin{Proposition}
\label{prop:converge:adelic:seq}
Let $X$ be an $e$-dimensional projective variety over $K$, and
let $L_1, \ldots, L_{e+1}$ be nef line bundles on $X$.
Let $\left\{ (\XX^{(i)}_n, \overline{\LL}^{(i)}_n ) \right\}$ be
an adelic sequence of $(X, L_i)$ for each $1 \leq i \leq e+1$.
Let $\ZZZ_n$ be a projective arithmetic variety over $B$
such that there are birational morphisms
$\mu_n^{(i)} : \ZZZ_n \to \XX^{(i)}_n$ \rom{(}$i=1, \ldots, e+1$\rom{)}.
Then, the limit
\[
\adeg\left( \acherncl_1({\mu_n^{(1)}}^*(\overline{\LL}^{(1)}_n)) \cdots
\acherncl_1({\mu_n^{(e+1)}}^*(\overline{\LL}^{(e+1)}_n)) \cdot
\acherncl_1(\pi^*_{\ZZZ_n}(\overline{H}_1)) \cdots 
\acherncl_1(\pi^*_{\ZZZ_n}(\overline{H}_d))
\right)
\]
as $n \to \infty$ exists, where $\pi_{\ZZZ_n} : \ZZZ_n \to B$ 
is the natural morphism.
Moreover, if $\left\{ (\YY^{(i)}_n, \overline{\MM}^{(i)}_n) \right\}$ is
another adelic sequence of $(X, L_i)$ for each $1 \leq i \leq e+1$, and
$\left\{ (\XX^{(i)}_n, \overline{\LL}^{(i)}_n) \right\}$ is equivalent to
$\left\{ (\YY^{(i)}_n, \overline{\MM}^{(i)}_n) \right\}$ for each $i$,
then the limit by $\left\{ (\XX^{(i)}_n, \overline{\LL}^{(i)}_n ) \right\}$
coincides with the limit by 
$\left\{ (\YY^{(i)}_n, \overline{\MM}^{(i)}_n) \right\}$.
\end{Proposition}

\Proof
Let $\ZZZ_{n,m}$ be a projective arithmetic variety over $B$ such that
there are birational morphisms $\ZZZ_{n,m} \to \ZZZ_n$ and
$\ZZZ_{n,m} \to \ZZZ_m$.
By abuse of notation, we denote birational morphisms 
$\ZZZ_{n, m} \to \XX^{(i)}_n$ and
$\ZZZ_{n, m} \to \XX^{(j)}_m$ by $\mu_n^{(i)}$ and $\mu_m^{(j)}$ respectively.
First of all, we can see
\begin{multline*}
\acherncl_1({\mu_n^{(1)}}^*(\overline{\LL}^{(1)}_n)) \cdots
\acherncl_1({\mu_n^{(e+1)}}^*(\overline{\LL}^{(e+1)}_n)) -
\acherncl_1({\mu_m^{(1)}}^*(\overline{\LL}^{(1)}_m)) \cdots
\acherncl_1({\mu_m^{(e+1)}}^*(\overline{\LL}^{(e+1)}_m)) \\
=
\sum_{i=1}^{e+1} \acherncl_1({\mu_n^{(1)}}^*(\overline{\LL}^{(1)}_n)) \cdots
\left( \acherncl_1({\mu_n^{(i)}}^*(\overline{\LL}^{(i)}_n)) -
\acherncl_1({\mu_m^{(i)}}^*(\overline{\LL}^{(i)}_m)) \right)
\cdots \acherncl_1({\mu_m^{(e+1)}}^*(\overline{\LL}^{(e+1)}_m))
\end{multline*}
Therefore, it is sufficient to show that, for any positive $\epsilon$,
there is a positive integer $N$ such that if $n, m \geq N$, then
\[
\adeg\left( \Delta_{n,m,i} \cdot
\acherncl_1(\pi^*_{\ZZZ_{n,m}}(\overline{H}_1)) \cdots 
\acherncl_1(\pi^*_{\ZZZ_{n,m}}(\overline{H}_d))
\right) \leq \epsilon,
\]
where
$\Delta_{n,m,i} = \acherncl_1({\mu_n^{(1)}}^*(\overline{\LL}^{(1)}_n))
\cdots
\left( \acherncl_1({\mu_n^{(i)}}^*(\overline{\LL}^{(i)}_n)) -
\acherncl_1({\mu_m^{(i)}}^*(\overline{\LL}^{(i)}_m)) \right)
\cdots \acherncl_1({\mu_m^{(e+1)}}^*(\overline{\LL}^{(e+1)}_m))$.
By the definition of adelic sequences,
there are a projective arithmetic variety $\XX_{n, m}$
over $B$, a birational morphism $\nu_{n, m} : \XX_{n,m} \to \ZZZ_{n,m}$,
and a nef $C^{\infty}$-hermitian $\QQ$-line bundle $\overline{D}_{n, m}$ on $B$
such that
\[
 -\pi_{\XX_{n,m}}^*(\overline{D}_{n, m}) \precsim
 \nu_{n,m}^*\left(
{\mu_n^{(i)}}^*(\overline{\LL}^{(i)}_n) - 
{\mu_m^{(i)}}^*(\overline{\LL}^{(i)}_m)
 \right)
  \precsim
 \pi_{\XX_{n,m}}^*(\overline{D}_{n, m}).
\]
Here, since $\overline{\LL}_n^{(i)}$'s are nef
with respect to $\XX_n^{(i)} \to B$ and $\overline{H}_j$'s are nef, 
by using Lemma~\ref{lem:non:negative:pi:nef} together with the projection
formula, we can see
\begin{multline*}
\left\vert
\adeg\left( \Delta_{n,m,i} \cdot
\acherncl_1(\pi^*_{\ZZZ_{n,m}}(\overline{H}_1)) \cdots 
\acherncl_1(\pi^*_{\ZZZ_{n,m}}(\overline{H}_d))
\right) \right\vert \\
\leq \deg(L_1 \cdots L_{i-1} \cdot L_{i+1} \cdots L_{e+1})
\left\vert \adeg\left( \acherncl_1(\overline{D}_{n,m}) \cdot
\acherncl_1(\overline{H}_1) \cdots \acherncl_1(\overline{H}_d)
\right) \right\vert.
\end{multline*}
Thus we get the first assertion.
The second one is obvious by the definition of equivalence.
\QED

\begin{Definition}[Adelic intersection number]
Let $\overline{L}_1, \ldots, \overline{L}_{e+1}$ be adelic line bundles on $X$.
Then, by the above proposition, 
the limit of intersection numbers does not depend on
the choice of adelic sequences representing each $\overline{L}_i$.
Thus, we may define the adelic intersection number
$\langle \overline{L}_1 \cdots \overline{L}_{e+1}
\rangle_{\overline{B}}$ to be the limit in 
Proposition~\ref{prop:converge:adelic:seq}.
\end{Definition}

Here let us consider the following two propositions.
The second  proposition is a property concerning the specialization
of adelic intersection number.

\begin{Proposition}
\label{prop:elem:prop:adelic:int}
Let $\overline{L}_1, \ldots, \overline{L}_{e+1}$ be
adelic line bundles on $X$.
Then, we have the following.
\begin{enumerate}
\renewcommand{\labelenumi}{(\arabic{enumi})}
\item
If $\overline{L}_1, \ldots, \overline{L}_{e+1}$ are nef, then
$\langle \overline{L}_1 \cdots \overline{L}_{e+1} \rangle_{\overline{B}} 
\geq 0$.

\item
Let $\overline{H}'_1, \cdots, \overline{H}'_d$ be nef $C^{\infty}$-hermitian line bundles
on $B$ with $\overline{H}'_i \succsim \overline{H}_i$ for all $i$.
If $\overline{L}_1, \ldots, \overline{L}_{e+1}$ are nef, then
\[
\langle \overline{L}_1 \cdots \overline{L}_{e+1}
\rangle_{(B; \overline{H}'_1, \cdots, \overline{H}'_d)}
\geq \langle \overline{L}_1 \cdots \overline{L}_{e+1}
\rangle_{(B; \overline{H}_1, \cdots, \overline{H}_d)}.
\]

\item
Let $g : Y \to X$ be a generically finite morphism of
projective varieties over $K$.
Then,
\[
\langle g^*(\overline{L}_1) \cdots g^*(\overline{L}_{e+1})
\rangle_{\overline{B}} = \deg(g)
\langle \overline{L}_1 \cdots \overline{L}_{e+1}
\rangle_{\overline{B}}
\]
\end{enumerate}
\end{Proposition}

\Proof
(1) is a consequence of (2) of Proposition~\ref{prop:nef:plus:ample}.
(2) follows from (4) of Proposition~\ref{prop:nef:plus:ample}.
(3) is a consequence of the projection formula.
\QED

\begin{Proposition}
\label{prop:inq:specialization}
Let $\{ (\XX_n, \overline{\LL}_n) \}$
be an adelic sequence of $(X, L)$ such that
$\overline{\LL}_n$ is nef for every $n$, and
let $\overline{L}$ be a nef adelic line bundle on $X$
given by the adelic sequence $\{ (\XX_n, \overline{\LL}_n) \}$.
Let $U$ be a common base of the adelic sequence 
$\{ (\XX_n, \overline{\LL}_n) \}$
\rom{(}cf. 
the definition
of adelic sequences in \rom{\S\S\ref{subsec:adelic:sequence}}\rom{)}.
Let $\gamma$ be a point of codimension one in $U_{\QQ}$ such that
$\XX_U$  is flat over $\gamma$ and
the fiber $X_{\gamma}$ of $\XX_U \to U$ over $\gamma$ is integral.
Then, $X_{\gamma}$ is a projective variety over the residue field 
$\kappa(\gamma)$ at $\gamma$, and
$L_{\gamma} = \rest{\LL}{X_{\gamma}}$ is a line bundle on $X_{\gamma}$.
Let $\Gamma$ be the Zariski closure of $\{ \gamma \}$ in $B$, and
$\ZZZ_n$ the Zariski closure of $X_{\gamma}$ in $\XX_n$.
If $\overline{H}_d$ is big, then we have the following.
\begin{enumerate}
\renewcommand{\labelenumi}{(\arabic{enumi})}
\item
$\left\{ (\ZZZ_n, \rest{\overline{\LL}_n}{\ZZZ_n}) \right\}$ 
is an adelic sequence of
$(X_{\gamma}, L_{\gamma})$ 
with respect to $\left(\Gamma; \rest{\overline{H}_1}{\Gamma}, \ldots,
\rest{\overline{H}_{d-1}}{\Gamma}\right)$.

\item
If we denote by $\overline{L}_{\gamma}$ the adelic line bundle arising from
the adelic sequence 
$\left\{ (\ZZZ_n, \rest{\overline{\LL}_n}{\ZZZ_n}) \right\}$,
then $\left\langle \overline{L}^{\cdot \dim X + 1} 
\right\rangle_{(B; \overline{H}_1, \ldots,\overline{H}_d)} = 0$
implies
$\left\langle  \overline{L}_{\gamma}^{\cdot \dim X_{\gamma} + 1}
\right\rangle_{\left(\Gamma; \rest{\overline{H}_1}{\Gamma}, \ldots, 
\rest{\overline{H}_{d-1}}{\Gamma}\right)} = 0$.
\end{enumerate}
\end{Proposition}

\Proof
First of all, by using Lemma~\ref{lem:exist:small:sec:zero},
we fix a positive integer $N$ and
a non-zero section $s \in H^0(B, H_d^{\otimes N})$ with
$s(\gamma) = 0$ and $\Vert s \Vert_{\sup} \leq 1$.
Then, $\zero(s) = \Gamma + \Sigma$ for some effective divisor $\Sigma$.

\medskip
(1) To prove (1), it is sufficient to show that
\[
\lim_{n,m \to \infty}
\adeg\left( \acherncl_1(\rest{\overline{D}_{n, m}}{\Gamma})
\cdot \acherncl_1(\rest{\overline{H}_1}{\Gamma}) \cdots 
\acherncl_1(\rest{\overline{H}_{d-1}}{\Gamma})
\right) = 0,
\]
where $\overline{D}_{n,m}$ is a nef $C^{\infty}$-hermitian 
$\QQ$-line bundle on $B$
appeared in the definition of adelic sequences 
(cf. \S\S\ref{subsec:adelic:sequence}).
First of all,
\begin{multline*}
N \adeg\left( \acherncl_1(\overline{D}_{n, m})
\cdot \acherncl_1(\overline{H}_1) \cdots 
\acherncl_1(\overline{H}_{d}) \right) \\
=
\adeg\left( \acherncl_1(\rest{\overline{D}_{n, m}}{\Gamma})
\cdot \acherncl_1(\rest{\overline{H}_1}{\Gamma}) \cdots 
\acherncl_1(\rest{\overline{H}_{d-1}}{\Gamma})
\right) \\
\qquad\qquad + \adeg\left( \acherncl_1(\rest{\overline{D}_{n, m}}{\Sigma})
\cdot \acherncl_1(\rest{\overline{H}_1}{\Sigma}) \cdots 
\acherncl_1(\rest{\overline{H}_{d-1}}{\Sigma})
\right) \\
+ \int_{B(\CC)}  -\log(\Vert s \Vert) 
c_1(\overline{D}_{n,m}) \wedge c_1(\overline{H}_1)
\wedge \cdots \wedge c_1(\overline{H}_{d-1}).
\end{multline*}
Here every term is non-negative. Thus, we can see that
\[
\lim_{n,m \to \infty}
\adeg\left( \acherncl_1(\rest{\overline{D}_{n, m}}{\Gamma})
\cdot \acherncl_1(\rest{\overline{H}_1}{\Gamma}) \cdots 
\acherncl_1(\rest{\overline{H}_{d-1}}{\Gamma})
\right) = 0.
\]

\medskip
(2)
We can set $\zero(\pi^{*}_{\XX_n}(s)) = \ZZZ_n + \Delta_n$ 
for some effective divisor $\Delta_n$.
Therefore,
\begin{multline*}
N \adeg\left(\acherncl_1(\overline{\LL}_n)^{\cdot e+1} \cdot 
\acherncl_1(\overline{H}_1) \cdots
\acherncl_1(\overline{H}_d) \right) \\
=
\adeg\left(\acherncl_1(\rest{\overline{\LL}_n}{\ZZZ_n})^{\cdot e+1} 
\cdot \acherncl_1(\rest{\pi_{\XX_n}^*(\overline{H}_1)}{\ZZZ_n}) \cdots
\acherncl_1(\rest{\pi_{\XX_n}^*(\overline{H}_{d-1})}{\ZZZ_n}) \right) \\
\qquad + \adeg\left(\acherncl_1(\rest{\overline{\LL}_n}{\Delta_n})^{\cdot e+1} 
\cdot \acherncl_1(\rest{\pi_{\XX_n}^*(\overline{H}_1)}{\Delta_n}) \cdots
\acherncl_1(\rest{\pi_{\XX_n}^*(\overline{H}_{d-1})}{\Delta_n}) \right) \\
\quad + \int_{\XX_n(\CC)} -\log(\Vert \pi_{\XX_n}^*(s) \Vert) 
c_1(\overline{\LL}_n)^{\wedge e+1} \wedge c_1(\pi_{\XX_n}^*(\overline{H}_1))
\wedge \cdots \wedge c_1(\pi_{\XX_n}^*(\overline{H}_{d-1})).
\end{multline*}
Since the last two terms of the above equation are non-negative, we have
\begin{multline*}
N \adeg\left(\acherncl_1(\overline{\LL}_n)^{\cdot e+1} \cdot 
\acherncl_1(\overline{H}_1) \cdots
\acherncl_1(\overline{H}_d) \right) \\
\geq
\adeg\left(\acherncl_1(\rest{\overline{\LL}_n}{\ZZZ_n})^{\cdot e+1} 
\cdot \acherncl_1(\rest{\pi_{\XX_n}^*(\overline{H}_1)}{\ZZZ_n}) \cdots
\acherncl_1(\rest{\pi_{\XX_n}^*(\overline{H}_{d-1})}{\ZZZ_n}) \right).
\end{multline*}
Thus, taking $n \to \infty$,
\[
N \left\langle \overline{L}^{\cdot \dim X + 1} 
\right\rangle_{(B; \overline{H}_1, \ldots,\overline{H}_d)}
\geq
\left\langle  \overline{L}_{\gamma}^{\cdot \dim X_{\gamma} + 1}
\right\rangle_{\left(\Gamma; \rest{\overline{H}_1}{\Gamma}, \ldots, 
\rest{\overline{H}_{d-1}}{\Gamma}\right)}.
\]
Therefore, we get (2).
\QED

\setcounter{Theorem}{0}
\subsection{Adelic height}
Let $K$ be a finitely generated field over $\QQ$
with $d = \trdeg_{\QQ}(K)$, 
and $\overline{B} = (B; \overline{H}_1, \ldots, \overline{H}_d)$
a polarization of $K$.
Let $X$ be a projective variety over $K$, and
$L$ an ample line bundle on $X$.

Let $\overline{L}$ be an adelic line bundle given by
an adelic sequence $\{ (\XX_n, \overline{\LL}_n) \}$.
Let $K'$ be a finite extension of $K$, $B'$ 
the normalization of $B$ in $K'$, and
let $\rho : B' \to B$ be the induced morphism.
Let $\XX'_n$ be the main component of $\XX_n \times_{B} B'$. 
We set the induced morphisms
as follows.
\[
\begin{CD}
\XX_n @<{\tau_n}<< \XX'_n \\
@V{\pi_n}VV @VV{\pi'_n}V \\
B     @<{\rho}<< B'
\end{CD}
\]
Then, $\{ \XX'_n, \tau_n^*(\overline{\LL}_n) \}$ 
is an adelic sequence of $(X_{K'}, L_{K'})$.
We denote by $\overline{L}_{K'}$ the adelic line bundle induced by
$\{ \XX'_n, \tau_n^*(\overline{\LL}_n) \}$.
With this notation, if $\overline{L}_1, \ldots, \overline{L}_{e+1}$
are adelic line bundles on $X$, then we can see
\renewcommand{\theequation}{\arabic{section}.\arabic{subsection}.\arabic{Theorem}}
\addtocounter{Theorem}{1}
\begin{equation}
\label{eqn:adelic:int:finite:extension}
\langle (\overline{L}_1)_{K'} \cdots (\overline{L}_{e+1})_{K'}
\rangle_{\overline{B}_{K'}} =
[K' : K]\langle \overline{L}_1 \cdots \overline{L}_{e+1}
\rangle_{\overline{B}},
\end{equation}
\renewcommand{\theequation}{\arabic{section}.\arabic{subsection}.\arabic{Theorem}.\arabic{Claim}}
by virtue of the projection formula,
where $\overline{B}_{K'} = 
(B'; \rho^*(\overline{H}_1), \cdots, \rho^*(\overline{H}_d))$.

\medskip
Let $Y$ be a subvariety of $X_{\overline{K}}$.
We assume that $Y$ is defined over $K'$.
Let $\YY_n'$ be the closure of $Y$ in $\XX'_n$.
Then, $\{ \YY'_n, \rest{\tau_n^*(\overline{\LL}_n)}{\YY'_n} \}$
is an adelic sequence of $(Y, \rest{L_{K'}}{Y})$.
We denote by $\rest{\overline{L}_{K'}}{Y}$ the adelic line bundle
given by $\{ \YY'_n, \rest{\tau_n^*(\overline{\LL}_n)}{\YY'_n} \}$.
We define the height of $Y$ with respect to $\overline{L}$ to be
\[
h_{\overline{L}}^{\overline{B}}(Y) = \frac{\left\langle 
\left(\rest{\overline{L}_{K'}}{Y}\right)^{\cdot \dim Y + 1} 
\right\rangle_{\overline{B}'}}%
{[K':K](\dim Y +1)\deg\left(\rest{L_{K'}}{Y}^{\dim Y}\right)}.
\]
Note that by virtue of \eqref{eqn:adelic:int:finite:extension},
the above does not depend on the choice of $K'$.
We call $h_{\overline{L}}^{\overline{B}}(Y)$ the 
{\em adelic height}
of $Y$ with respect to $\overline{L}$ and $\overline{B}$.

\begin{Proposition}
\label{prop:height:zero:period}
Let $X$ be a projective variety over $K$, and $L$ an ample line bundle on $X$.
We assume that there is a surjective morphism $f : X \to X$
and an integer $d \geq 2$
with $L^{\otimes d} \simeq f^*(L)$.
Let $\overline{L}^f$ be the adelic line bundle with the $f$-adelic structure.
Then, we have the following.
\begin{enumerate}
\renewcommand{\labelenumi}{(\arabic{enumi})}
\item
$h^{\overline{B}}_{\overline{L}^f}(Y) \geq 0$
for all subvarieties $Y$ of $X_{\overline{K}}$.

\item
For a $C^{\infty}$-model $(\XX, \overline{\LL})$ of $(X, L)$
with $\overline{\LL}$ being nef with respect to $\XX \to B$,
there is a constant $C$ such that
\[
\left\vert h^{\overline{B}}_{\overline{L}^f}(Y) -
h^{\overline{B}}_{(\XX, \overline{\LL})}(Y) \right\vert \leq C
\]
for any subvarieties $Y$ of $X_{\overline{K}}$.

\item
$h^{\overline{B}}_{\overline{L}^f}(f(Y)) = 
d h^{\overline{B}}_{\overline{L}^f}(Y)$ for
any subvarieties $Y$ of $X_{\overline{K}}$.
\end{enumerate}
Moreover, $h^{\overline{B}}_{\overline{L}^f}$ is characterized by the above properties
\rom{(1)}, \rom{(2)} and \rom{(3)}.
\end{Proposition}

\Proof
(1) Since $\overline{L}^f$ is nef by Proposition~\ref{prop:f:adelic:sequence}, 
(1) is a consequence of (1) of Proposition~\ref{prop:elem:prop:adelic:int}.

\medskip
(2) We choose a Zariski open set $U$ of $B$ such that
$f$ extends to $f_U: \XX_U \to \XX_U$ and
$\LL_U^{\otimes d} = f^*(\LL_U)$
in $\Pic(\XX_U) \otimes \QQ$.
Let $\XX_n$ be the normalization of 
$\XX_U \overset{f_U^n}{\longrightarrow} \XX_U \to \XX$, and
$f_n : \XX_n \to \XX$ the induced morphism.
We denote $f_n^*(\overline{\LL})^{\otimes d^{-n}}$ 
by $\overline{\LL}_n$.
Then, as in proof of (1) of Proposition~\ref{prop:f:adelic:sequence}, there
are a projective arithmetic variety $\ZZZ_n$ over $B$,
birational morphisms $\mu_n : \ZZZ_n \to \XX_n$ and $\nu_n : \ZZZ_n \to \XX$
(which are the identity map over $U$), and
an ample $C^{\infty}$-hermitian line bundle
$\overline{D}$ on $B$ such that
\[
 - d_n \pi_{\ZZZ_n}^*(\overline{D}) 
\precsim_{\pi_{\ZZZ_n}^{-1}(U)}
\mu_n^*(\overline{\LL}_n) - \nu_n^*(\overline{\LL})
\precsim_{\pi_{\ZZZ_n}^{-1}(U)}
 d_n \pi_{\ZZZ_n}^*(\overline{D}),
\]
where $d_n = \sum_{j=0}^{n-1} d^{-j}$.

Let $Y$ be a subvariety of $X_{\overline{K}}$.
We assume that $Y$ is defined over a finite extension field $K'$ of $K$.
Let $B'$ be
the normalization of $B$ in $K'$, and
let $\rho : B' \to B$ be the induced morphism.
We denote by $\XX'$, $\XX'_n$ and $\ZZZ'_n$ the
main components of $\XX \times_{B} B'$, $\XX_n \times_{B} B'$ and
$\ZZZ_n \times_{B} B'$ respectively.
We set the induced morphisms as follows.
\[
\begin{CD}
\XX @<{\tau}<< \XX' \\
@V{\pi_{\XX}}VV @VV{\pi_{\XX'}}V \\
B     @<{\rho}<< B'
\end{CD}
\qquad\qquad
\begin{CD}
\XX_n @<{\tau_n}<< \XX'_n \\
@V{\pi_{\XX_n}}VV @VV{\pi_{\XX'_n}}V \\
B     @<{\rho}<< B'
\end{CD}
\qquad\qquad
\begin{CD}
\ZZZ_n @<<< \ZZZ'_n \\
@V{\pi_{\ZZZ_n}}VV @VV{\pi_{\ZZZ'_n}}V \\
B     @<{\rho}<< B'
\end{CD}
\]
We also have the induced morphisms
$\mu'_n : \ZZZ'_n \to \XX'_n$ and $\nu'_n : \ZZZ'_n \to \XX'$.
Then,
\[
 - d_n \pi_{\ZZZ'_n}^*(\rho^*\overline{D}) 
\precsim_{\pi_{\ZZZ_n}^{-1}(U)}
{\mu'_n}^*(\tau_n^* \overline{\LL}_n) - {\nu'_n}^*(\tau^* \overline{\LL})
\precsim_{\pi_{\ZZZ_n}^{-1}(U)}
 d_n \pi_{\ZZZ'_n}^*(\rho^*\overline{D}),
\]
On the other hand, since
\begin{multline*}
\acherncl_1({\mu'_n}^*(\tau_n^* \overline{\LL}_n))^{\dim Y + 1} - 
\acherncl_1({\nu'_n}^*(\tau^* \overline{\LL}))^{\dim Y + 1} \\
= \sum_{i=1}^{\dim Y + 1} \acherncl_1({\mu'_n}^*(\tau_n^* \overline{\LL}_n))^{i-1}
\cdot \left(  \acherncl_1({\mu'_n}^*(\tau_n^* \overline{\LL}_n)) - 
\acherncl_1({\nu'_n}^*(\tau^* \overline{\LL}))
\right) \acherncl_1({\nu'_n}^*(\tau^* \overline{\LL}))^{\dim Y - i + 1},
\end{multline*}
by using Lemma~\ref{lem:non:negative:pi:nef}, we have
\begin{multline*}
\left\vert
\adeg \left( \acherncl_1(\rest{\tau_n^* \overline{\LL}_n}{\YY_n})^{\dim Y + 1} \cdot
\acherncl_1(\pi_{\YY_n}^*\rho^* \overline{H}_1) \cdots 
\acherncl_1(\pi_{\YY_n}^*\rho^* \overline{H}_d)
\right) \right.\\
\left. -\adeg \left( \acherncl_1(\rest{\tau^* \overline{\LL}}{\YY})^{\dim Y + 1} \cdot
\acherncl_1(\pi_{\YY}^*\rho^* \overline{H}_1) \cdots 
\acherncl_1(\pi_{\YY}^*\rho^* \overline{H}_d)
\right) \right\vert \\
\leq d_n [K':K](\dim Y + 1)
\deg(\rest{L}{Y}^{\dim Y}) \adeg\left( \acherncl_1(\overline{D}) \cdot
\acherncl_1(\overline{H}_1) \cdots \acherncl_1(\overline{H}_d) \right),
\end{multline*}
where $\YY$ and $\YY_n$ are the Zariski closures of $Y$ in $\XX'$ and $\XX'_n$
respectively.
Thus we get (2).

\medskip
(3) Clearly, we may assume $Y$ is defined over $K$.
Let $(\XX, \overline{\LL})$ be a $C^{\infty}$ model of $(X, L)$.
Let us consider a sequence of morphisms of projective arithmetic varieties over $B$:
\[
\XX = \XX_0 \overset{f_1}{\longleftarrow} \XX_1 
\overset{f_2}{\longleftarrow} \cdots \overset{f_{n-1}}{\longleftarrow}
\XX_{n-1} \overset{f_{n}}{\longleftarrow} \XX_n 
\overset{f_{n+1}}{\longleftarrow} \XX_{n+1} \overset{f_{n+2}}{\longleftarrow}
\cdots
\]
such that $X$ is the generic fiber of $\XX_n \to B$ for every $n$, and
that
$f_n : \XX_n \to \XX_{n-1}$ is an extension of $f$
for each $n$.
Let $\YY_n$ be the Zariski closure of $Y$ in $\XX_n$.
Then, $f_{n+1}(\YY_{n+1})$ is the Zariski closure of $f(Y)$ in
$\XX_n$.
By the definition of the height,
\addtocounter{Claim}{1}
\begin{multline}
\label{prop:height:zero:period:eqn:1}
h^{\overline{B}}_{\overline{L}^f}(Y) \\
= \lim_{n\to\infty} 
\frac{ \adeg\left( 
\acherncl_1(f_{n+1}^* f_n^* \cdots f_1^*(\overline{\LL}))^{\cdot \dim Y +1} 
\cdot
\acherncl_1(f_{n+1}^* \pi_{\XX_n}^*(\overline{H}_1)) \cdots
\acherncl_1(f_{n+1}^* \pi_{\XX_n}^*(\overline{H}_d)) \cdot (\YY_{n+1},0) 
\right)}
{(\dim Y +1)\deg(\rest{L}{Y}^{\dim Y}) d^{(n+1)(\dim Y+1)}}.
\end{multline}
On the other hand, by the projection formula,
\addtocounter{Claim}{1}
\begin{multline}
\label{prop:height:zero:period:eqn:2}
\adeg\left( 
\acherncl_1(f_{n+1}^* f_n^* \cdots f_1^*(\overline{\LL}))^{\cdot \dim Y+1} 
\cdot
\acherncl_1(f_{n+1}^* \pi_{\XX_n}^*(\overline{H}_1)) \cdots
\acherncl_1(f_{n+1}^* \pi_{\XX_n}^*(\overline{H}_d)) \cdot (\YY_{n+1},0) 
\right) \\
= \deg (\rest{f}{Y}) \adeg\left( 
\acherncl_1(f_n^* \cdots f_1^*(\overline{\LL}))^{\cdot \dim Y+1} \cdot
\acherncl_1(\pi_{\XX_n}^*(\overline{H}_1)) \cdots
\acherncl_1(\pi_{\XX_n}^*(\overline{H}_d)) \cdot (f_{n+1}(\YY_{n+1}), 0) 
\right).
\end{multline}
Here, since $L^{\otimes d} \simeq f^*(L)$, 
we have $\rest{L^{\otimes d}}{Y} \simeq 
(\rest{f}{Y})^*\left( \rest{L}{f(Y)} \right)$, which implies
\addtocounter{Claim}{1}
\begin{equation}
\label{prop:height:zero:period:eqn:3}
d^{\dim Y} \deg(\rest{L}{Y}^{\dim Y}) = 
\deg(\rest{f}{Y}) \deg(\rest{L}{f(Y)}^{\dim Y}).
\end{equation}
Moreover,
\addtocounter{Claim}{1}
\begin{multline}
\label{prop:height:zero:period:eqn:4}
h^{\overline{B}}_{\overline{L}^f}(f(Y)) \\
= \lim_{n\to\infty} 
\frac{\adeg\left( 
\acherncl_1(f_n^* \cdots f_1^*(\overline{\LL}))^{\cdot \dim Y+1} \cdot
\acherncl_1(\pi_{\XX_n}^*(\overline{H}_1)) \cdots
\acherncl_1(\pi_{\XX_n}^*(\overline{H}_d)) \cdot (f_{n+1}(\YY_{n+1}), 0) \right)}
{(\dim Y +1)\deg(\rest{L}{f(Y)}^{\dim Y}) d^{n(\dim Y+1)}}.
\end{multline}
Therefore, by \eqref{prop:height:zero:period:eqn:1},
\eqref{prop:height:zero:period:eqn:2},
\eqref{prop:height:zero:period:eqn:3}, and
\eqref{prop:height:zero:period:eqn:4}, we obtain
\[
h^{\overline{B}}_{\overline{L}^f}(f(Y)) = 
d h^{\overline{B}}_{\overline{L}^f}(Y).
\]

\medskip
Finally, the last assertion is obvious. For, by (2) and (3),
we can see
\[
h_{\overline{L}^f}^{\overline{B}}(Y) =
\lim_{n \to \infty} \frac{h_{(\XX, \overline{L})}^{\overline{B}}(f^{\cdot n}(Y))}{d^n}.
\]
\QED

\renewcommand{\theTheorem}{\arabic{section}.\arabic{Theorem}}
\renewcommand{\theClaim}{\arabic{section}.\arabic{Theorem}.\arabic{Claim}}
\renewcommand{\theequation}{\arabic{section}.\arabic{Theorem}.\arabic{Claim}}
\section{The canonical height of subvarieties of
an abelian variety over finitely generated fields}
Let $K$ be a finitely generated field over $\QQ$ with $d = \trdeg_{\QQ}(K)$, 
and $\overline{B} = (B; \overline{H}_1, \ldots, \overline{H}_d)$
a polarization of $K$.
Let $A$ be an abelian variety over $K$, and 
$L$ a symmetric ample line bundle on $A$.
Since $[2]^*(L) \simeq L^{\otimes 4}$, we have 
an adelic line bundle $\overline{L}^{[2]}$ with the $[2]$-adelic structure.
Let $f : A \to A$ be an endomorphism with $f^*(L) \simeq L^{\otimes d}$
for some $d \geq 2$.
Then, since $f \cdot [2] = [2] \cdot f$, by
(2) of Proposition~\ref{prop:f:adelic:sequence}, 
$\overline{L}^{f} = \overline{L}^{[2]}$.
Thus, the adelic structure does not depend on the choice of
the endomorphism.
In this sense, we have the line bundle $\overline{L}^{can}$
with the canonical adelic structure.

Let $X$ be a subvariety of $A_{\overline{K}}$.
We denote by $\hat{h}^{\overline{B}}_L(X)$
the adelic height $h_{\overline{L}^{can}}^{\overline{B}}(X)$ of $X$ with
respect to the line bundle $\overline{L}^{can}$ with
the canonical adelic structure.
Then, by Proposition~\ref{prop:height:zero:period},
we can see the following:
\begin{enumerate}
\renewcommand{\labelenumi}{(\alph{enumi})}
\item
$\hat{h}_L^{\overline{B}}(X) \geq 0$
for all subvarieties $X$ of $A_{\overline{K}}$.

\item
For a $C^{\infty}$-model $(\mathcal{A}, \overline{\LL})$ of $(A, L)$
with $\overline{\LL}$ being nef with respect to $\mathcal{A} \to B$,
there is a constant $C$ such that
\[
 \left\vert \hat{h}_L^{\overline{B}}(X) - 
        h^{\overline{B}}_{(\mathcal{A}, \overline{\LL})}(X) \right\vert \leq C
\]
for all subvarieties $X$ of $A_{\overline{K}}$.

\item
$\hat{h}_L^{\overline{B}}([N](X)) = N^2 \hat{h}_L^{\overline{B}}(X)$
for all subvarieties $X$ of $A_{\overline{K}}$ and
all non-zero integers $N$.
\end{enumerate}

\medskip
The purpose of this section is to prove the following theorem.

\begin{Theorem}
\label{thm:height:zero}
Let $A$ be an abelian variety over $K$, and
$L$ a symmetric ample line bundle on $A$.
Let $X$ be a subvariety of $A_{\overline{K}}$.
If the polarization $\overline{B}$ is big, then the following are
equivalent.
\begin{enumerate}
\renewcommand{\labelenumi}{(\arabic{enumi})}
\item
$X$ is a translation of an abelian subvariety by a torsion point.

\item
The set 
$\{ x \in X(\overline{K}) \mid \hat{h}_L^{\overline{B}}(x)
   \leq \epsilon \}$
is Zariski dense in $X$ for every $\epsilon > 0$.

\item
$\hat{h}^{\overline{B}}_{L}(X) = 0$.
\end{enumerate}
\end{Theorem}

\Proof
Let us begin with the following two lemmas.

\begin{Lemma}
\label{lem:comp:height:quotient}
Let $A$ be an abelian subvariety over $K$, $C$ an abelian subvariety of $A$,
and $\rho : A \to A' = A/C$ the natural homomorphism.
Let $X$ be a subvariety of $A$ such that $X = \rho^{-1}(\rho(X))$.
Let $L$ and $L'$ be symmetric ample line bundles on $A$ and $A'$ respectively.
If $\hat{h}_{L}^{\overline{B}}(X) = 0$, then
$\hat{h}_{L'}^{\overline{B}}(Y) = 0$, where $Y = \rho(X)$.
\end{Lemma}

\Proof
Replacing $L$ by $L^{\otimes n}$ ($n > 0$), we may assume that
$L \otimes \rho^{*}(L')^{\otimes -1}$ is generated by global sections.
Let $(\mathcal{A}, \overline{\LL})$ and $(\mathcal{A}', \overline{\LL}')$
be $C^{\infty}$-models of $(A, L)$ and $(A', L')$ 
over $B$ with the following properties:
\begin{enumerate}
\renewcommand{\labelenumi}{(\arabic{enumi})}
\item
$\overline{\LL}$ and $\overline{\LL}'$ are nef and big.

\item
There is a morphism $\mathcal{A} \to \mathcal{A'}$ over $B$ as an extension of
$\rho : A \to A'$.
(By abuse of notation, the extension is also denoted by $\rho$.)
\end{enumerate}
Let $\pi : \mathcal{A} \to B$ be the canonical morphism.
Replacing $\LL$ by $\LL \otimes \pi^{*}(Q)$ for some ample line bundle $Q$
on $B$, we may assume that $\pi_*(\LL \otimes \rho^*(\LL')^{\otimes -1})$
is generated by global sections.
Thus, there are sections $s_1, \ldots, s_r$ of 
$H^0(\LL \otimes \rho^*(\LL')^{\otimes -1})$
such that $\{ s_1, \ldots, s_r \}$ generates 
$L \otimes \rho^{*}(L')^{\otimes -1}$ on $A$.
Moreover, replacing the metric of $\overline{\LL}$, we may assume that
$s_1, \ldots, s_r$ are small sections, i.e.,
$\Vert s_i \Vert_{\sup} < 1$ for all $i$.

Let $\mathcal{A}_n$ (resp. $\mathcal{A}'_n$) be the normalization of
$A \overset{[2^n]}{\longrightarrow} A \hookrightarrow \mathcal{A}$
(resp. $A' \overset{[2^n]}{\longrightarrow} A' \hookrightarrow \mathcal{A}'$).
Then, we have the following commutative diagram:
\[
\begin{CD}
\mathcal{A} @<{f_n}<<  \mathcal{A}_n \\
@V{\rho}VV             @VV{\rho_n}V           \\  
\mathcal{A}' @<<{f'_n}< \mathcal{A}'_n, \\
\end{CD}
\]
where $f_n$ and $f'_n$ are extension of $[2^n]$.
Here the adelic structure of $\overline{L}$ (resp. $\overline{L}'$)
is induced by 
$\{ 4^{-n} f_n^*(\overline{\LL}) \}$ 
(resp. $\{ 4^{-n} {f_n'}^*(\overline{\LL}') \}$).
Let $\XX_n$ (resp. $\YY_n$) be the Zariski closure of $X$ in $\mathcal{A}_n$
(resp. $Y$ in $\mathcal{A}'_n$).
Then, since $f_n^*(s_1), \ldots, f_n^*(s_r)$ generate 
$f_n^*(L \otimes \rho^*(L')^{\otimes -1})$
on $A$, we can find $f_n^*(s_i)$ such that $f_n^*(s_i) \not= 0$ on $\XX_n$.
This means that 
$\rest{f_n^*(\overline{\LL})}{\XX_n} \otimes 
\rest{\rho_n^*({f'_n}^*(\overline{\LL}'))^{\otimes -1}}{\XX_n}$
is effective.
Therefore, if we denote $\dim X$ and $\dim Y$ by
$e$ and $e'$ respectively, then, by virtue of (4) of
Proposition~\ref{prop:nef:plus:ample} together 
with the projection formula,
\begin{multline*}
\adeg\left(
\acherncl_1(\rest{f_n^*(\overline{\LL})}{\XX_n})^{\cdot e+1} \cdot 
\acherncl_1(\pi_{\XX_n}^*\overline{H}_1)
\cdots \acherncl_1(\pi_{\XX_n}^*\overline{H}_d) \right) \\
\geq \adeg\left(
\acherncl_1(\rest{\rho_n^*{f'_n}^*(\overline{\LL}')}{\XX_n})^{\cdot e'+1} \cdot 
\acherncl_1(\rest{f_n^*(\overline{\LL})}{\XX_n})^{\cdot e-e'} \cdot
\acherncl_1(\rho_n^*\pi_{\YY_n}^*\overline{H}_1)
\cdots \acherncl_1(\rho_n^*\pi_{\YY_n}^*\overline{H}_d) \right) \\
= 4^{n(e-e')}\deg(\rest{L}{C}^{e-e'}) \adeg\left(
\acherncl_1(\rest{{f'_n}^*(\overline{\LL}')}{\YY_n})^{\cdot e'+1}  \cdot
\acherncl_1(\pi_{\YY_n}^*\overline{H}_1)
\cdots \acherncl_1(\pi_{\YY_n}^*\overline{H}_d) \right).
\end{multline*}
Hence,
\[
\hat{h}^{\overline{B}}_{\overline{L}}(X) \geq 
\frac{(e'+1) \deg(\rest{L'}{Y}^{e'}) \deg(\rest{L}{C}^{e-e'})}
{(e+1)\deg(\rest{L}{X}^e)} \hat{h}^{\overline{B}}_{\overline{L}'}(Y).
\]
Thus we get our assertion.
\QED

\begin{Lemma}
\label{lem:trans:abelian:family}
Let $A$ and $S$ be algebraic varieties over a field of characteristic zero, and
let $f : A \to S$ be an abelian scheme. Let $X$ be a subvariety of $A$ 
such that
$\rest{f}{X} : X \to B$ is proper and flat.
Let $s$ be a point of $S$.
If $X_{\bar{s}}$ is a translation of an abelian subvariety of $A_{\bar{s}}$,
then there is a Zariski open set $U$ of $S$ such that
\rom{(1)} $s \in U$ and \rom{(2)} $X_{\bar{t}}$ is 
a translation of an abelian subvariety of $A_{\bar{t}}$
for all $t \in U$.
In particular, the geometric generic fiber $X_{\bar{\eta}}$ is
a translation of an abelian subvariety. 
\end{Lemma}

\Proof
Since $X_{\bar{s}}$ is smooth and $q(X_{\bar{s}}) = \dim (X/S)$,
there is a Zariski open set $U$ of $S$ such that
$s \in U$, $X_U$ is smooth over $U$, and that
$q(X_{\bar{t}}) \leq \dim(X/S)$ for all $t \in U$.
By Ueno's theorem (cf. \cite[Theorem~10.12]{Iitaka}),
$q(X_{\bar{t}}) \geq \dim(X/S)$ and the equality holds if and only if
$X_{\bar{t}}$ is a a translation of an abelian subvariety. 
Thus we get our lemma.
\QED

\bigskip
Let us start the proof of Theorem~\ref{thm:height:zero}.
First of all, we may assume that $X$ is defined over $K$.

``(1) $\Longrightarrow$ (2)'' is obvious.
``(2) $\Longrightarrow$ (1)'' is nothing more than
Bogomolov's conjecture solved in \cite{MoArht}.

\medskip
``(1) $\Longrightarrow$ (3)'':
We set $X = A' + x$, where $A'$ is an abelian subvariety of $A_{\overline{K}}$ and
$x$ is a torsion point.
Let $N$ be a positive integer with $N x = 0$ and $N \geq 2$.
Then, $[N](X) = A' = [N](A')$.
Thus, by Proposition~\ref{prop:height:zero:period},
\[
\hat{h}^{\overline{B}}_{L}(X) = (1/N^2)\hat{h}^{\overline{B}}_{L}([N](X)) =
(1/N^2) \hat{h}^{\overline{B}}_{L}([N](A')) = \hat{h}^{\overline{B}}_{L}(A').
\]
On the other hand,
\[
\hat{h}^{\overline{B}}_{L}(A') = \hat{h}^{\overline{B}}_{L}([N](A')) =
N^2\hat{h}^{\overline{B}}_{L}(A').
\]
Therefore, $\hat{h}^{\overline{B}}_{L}(X) = \hat{h}^{\overline{B}}_{L}(A') = 0$.

\medskip
``(3) $\Longrightarrow$ (1)'': 
Let $\overline{H}$ be an ample $C^{\infty}$-hermitian line bundle on $B$.
Then, there is a positive integer $n$ such that
$\overline{H}_i^{\otimes n} \otimes \overline{H}^{\otimes -1} \succsim 0$.
for all $i$.
Then, by using (4) of Proposition~\ref{prop:nef:plus:ample},
we can see that an adelic sequence with respect to
$(B; \overline{H}_1, \ldots, \overline{H}_d)$ 
is an adelic sequence with respect to
$(B; \overline{H}, \ldots, \overline{H})$, and that
\[
0 \leq \hat{h}^{(B; \overline{H}, \ldots, \overline{H})}_{L}(X) \leq
n^d \hat{h}^{(B; \overline{H}_1, \ldots, \overline{H}_d)}_{L}(X).
\]
Thus, we may assume that $\overline{H}_1, \ldots, \overline{H}_d$ are ample.
We prove the assertion ``(3) $\Longrightarrow$ (1)''
by induction on $d = \trdeg_{\QQ}(K)$.
If $d = 0$, then this was proved by Zhang \cite{ZhEqui}.
We assume $d > 0$. Then, by the above lemma together with
hypothesis of induction and Proposition~\ref{prop:inq:specialization}, 
$X$ is a translation of an abelian subvariety $C$.
Let us consider $\pi : A \to A' = A/C$.
Then, $\pi(X)$ is a point, say $P$.
Then, by Lemma~\ref{lem:comp:height:quotient},
$\hat{h}_{L'}^{\overline{B}}(P) = 0$ 
for a symmetric ample line bundle $L'$ on $A'$.
Thus, $P$ is a torsion point by \cite[Proposition~3.4.1]{MoArht}.
Therefore, we can see that $X$ is a translation of $C$ by a torsion point.
\QED

Let $X$ be a smooth projective curve of genus $g \geq 2$ over 
$K$. Let $J$ be the Jacobian of $X$ and $L_{\Theta}$ a line bundle
given by a symmetric theta
divisor $\Theta$ on $J$, i.e., $L_{\Theta} = \OO_J(\Theta)$.
Let $j : X \to J$ be a morphism given by
$j(x) = \omega_X - (2g-2)x$. Then, it is well known that
$j^*(L_{\Theta}) = \omega_X^{\otimes 2g(g-1)}$.
Let $\overline{L}_{\Theta}^{can}$ be 
the canonical adelic structure of $L_{\Theta}$. 
Thus, we have the adelic line bundle
$j^*(\overline{L}_{\Theta}^{can})$ on $X$. 
In terms of this, we can give
the canonical adelic structure on $\omega_X$. We denote this by
$\overline{\omega}^{a}_X$. Then, as a corollary of
Theorem~\ref{thm:height:zero} and (3) of
Proposition~\ref{prop:elem:prop:adelic:int},  we have the following.

\begin{Corollary}
\label{cor:positive:self:int:can}
If the polarization $\overline{B}$ is big,
then
$\left\langle 
\overline{\omega}^{a}_X \cdot \overline{\omega}^{a}_X
\right\rangle_{\overline{B}} > 0$.
\end{Corollary}

%%%
%% Appendix
%%%
%\renewcommand{\thesection}{Appendix \Alph{section}}
%\renewcommand{\theTheorem}{\Alph{section}.\arabic{Theorem}}
%\renewcommand{\theClaim}{\Alph{section}.\arabic{Theorem}.\arabic{Claim}}
%\renewcommand{\theequation}{\Alph{section}.\arabic{Theorem}.\arabic{Claim}}
% beginning from A
%\setcounter{section}{0}
%\section{Local description of adelic structure}

\end{document}